\documentclass[a4paper, 12pt]{article}
\usepackage[T2A]{fontenc}
\usepackage[english,russian]{babel}
\usepackage{amsmath}
\usepackage{amsfonts}
\DeclareMathAlphabet{\mathpzc}{OT1}{pzc}{m}{it}

\newcommand{\Ker}{\operatorname{Ker}}

\newcommand{\e}{\mathfrak e}
\renewcommand{\l}{\mathfrak l}

\newcommand{\card}{\operatorname{card}}

\newcommand{\supp}{\operatorname{supp}}
\newcommand{\supvrai}{\operatornamewithlimits{sup\,vrai}}
\newcommand{\N}{\mathbb N}
\newcommand{\Z}{\mathbb Z}
\newcommand{\R}{\mathbb R}
\newcommand{\Nu}{\mathcal N}
\newcommand{\D}{\mathcal D}
\newcommand{\J}{\mathcal J}
\newcommand{\I}{\mathcal I}
\newcommand{\Kappa}{\mathcal K}
\renewcommand{\L}{\mathcal L}

\allowdisplaybreaks
\begin{document}

\author{ С. Н. Кудрявцев }
\title{Бернштейновские поперечники неизотропных классов
Никольского -- Бесова}
\date{}
\maketitle
\begin{abstract}
В статье рассмотрены неизотропные пространства Никольского и Бесова с
нормами, в определении которых вместо модулей непрерывности известных порядков
производных функций по координатным направлениям используются
"$L_p$-усредненные" модули непрерывности функций соответствующих порядков по
тем же направлениям. Для единичных шаров таких пространств функций, заданных в
областях определенного типа, получена слабая асимптотика поведения
$n$-поперечников по~Бернштейну в $L_q$-пространствах.
\end{abstract}

Ключевые слова: неизотропные пространства Никольского -- Бесова,
поперечник по Бернштейну
\bigskip

\centerline{Введение}
\bigskip

Описание асимптотики поведения бернштейновского $n$-поперечника
различных классов функций конечной гладкости в пространствах с интегральной
нормой дано в ряде публикаций. Среди них отметим [1] -- [4] (см. также
приведенную там литературу). Однако автору не приходилось встречать статьи,
в которых изучается поведение $n$-поперечника по Бернштейну неизотропных
классов Никольского и Бесова.

   Настоящая работа распространяет исследования асимптотики
поведения $n$-попе\-речника по~Бернштейну на неизотропные классы
функций Никольского и Бесова в пространствах с интегральной нормой.

А именно, при $ d \in \N, \alpha \in \R_+^d, 1 \le p,\theta < \infty $ для области $ D
\subset \R^d $ вводятся в рассмотрение
пространства $ (B_{p,\theta}^\alpha)^\prime(D) ((H_p^\alpha)^\prime(D)) $
с нормами
$$
\| f \|_{(B_{p,\theta}^\alpha)^\prime(D)} = \max\biggl(\| f \|_{L_p(D)},
\max_{j =1,\ldots,d} \left(\int_0^\infty
t^{-1 -\theta \alpha_j} (\Omega_j^{\prime l_j}(f, t)_{L_p(D)})^{\theta}\,dt
\right)^{1/\theta}\biggr),
$$
$$
\| f \|_{(H_p^\alpha)^\prime(D)} = \max(\| f \|_{L_p(D)}, \max_{j =1, \ldots, d}
\sup_{t \in \R_+} t^{-\alpha_j} \Omega_j^{\prime l_j}(f,t)_{L_p(D)}),
$$
где
\begin{multline*}
\Omega_j^{\prime l_j}(f,t)_{L_p(D)} =
\biggl((2 t)^{-1} \int_{ |\xi| \le t}
\| \Delta_{\xi e_j}^{l_j} f\|_{L_p(D_{l_j \xi e_j})}^p d\xi\biggr)^{1 /p},
t \in \R_+, l_j = \\
\min \{m \in \N: \alpha_j < m \}, D_{l_j \xi e_j} \ \text{см. в п.} \  1.1., j =1,\ldots,d.
\end{multline*}
В случае ограниченных областей $D $ так называемого $ \alpha$-типа, при
определенных условиях установлены верхняя и нижняя оценки бернштейновского
$n$-поперечника единичного шара пространства $ (B_{p,\theta}^\alpha)^\prime(D)
((H_p^\alpha)^\prime(D)) $ в пространстве $ L_q(D), $ дающие слабую
асимптотику этой величины.
Показано, что при $ C = B((B_{p,\theta}^\alpha)^\prime(D)),
B((H_p^\alpha)^\prime(D)), X = L_q(D) $ для бернштейновского $n$-поперечника
имеет место соотношение
\begin{multline*}
b_n(C,X) \asymp
\begin{cases}
n^{-1 /(\alpha^{-1}, \e)}, \parbox[t]{7cm}{ при $ 1 \le q \le p \le 2$ или $ 1 \le q = p <
\infty$ или ( $ 1 \le p < q \le \infty$ и $ 1 -(\alpha^{-1}, \e) /p +
(\alpha^{-1}, \e) /q >0 $ );} \\
n^{-\!(1 /(\alpha^{-1}, \e) -1 /p +1/2)},  \parbox[t]{5.5cm}{ при $ 1 \le q \le 2 < p < \infty $
и $ 1 -(\alpha^{-1}, \e) /p > 0$;} \\
n^{-\!(1 /(\alpha^{-1}, \e) -1 /p +1 /q)},  \parbox[t]{5.5cm}{ при $ 2 \le q < p < \infty $
и $ 1 -(\alpha^{-1}, \e) /p +(\alpha^{-1}, \e) /q -(\alpha^{-1}, \e) /2 >0,$}
\end{cases}
\end{multline*}
где $ (\alpha^{-1}, \e) = \sum_{j =1}^d 1 / \alpha_j. $
При этом для вывода верхней и нижней  оценки рассматриваемой величины
использованы другие средства, чем те, что применялись в указанных выше
работах. Кроме того, потребовались дополнительные элементы к схемам вывода
верхней и нижней оценок изучаемого поперечника.

Работа состоит из введения и двух параграфов, в первом из
которых даны предварительные сведения и вспомогательные результаты, а
во втором проведена соответственно верхняя и нижняя
оценка интересующей нас величины.
\bigskip

\centerline{\S 1. Предварительные сведения и вспомогательные
утверждения}
\bigskip

1.1. В этом пункте вводятся обозначения, относящиеся к
функциональным пространствам и классам, рассматриваемым в настоящей работе,
а также приводятся некоторые факты, необходимые в дальнейшем.

Для $ d \in \N $ через $ \Z_+^d $ обозначим множество
$$
\Z_+^d = \{\lambda = (\lambda_1, \ldots, \lambda_d) \in \Z^d:
\lambda_j \ge0, j =1, \ldots, d\}.
$$
Обозначим также при $ d \in \N $ для $ l \in \Z_+^d $ через
$ \Z_+^d(l) $ множество
\begin{equation*}
\Z_+^d(l) = \{ \lambda \in \Z_+^d: \lambda_j \le l_j, j =1, \ldots,d\}.
\end{equation*}

Напомним, что при $ n \in \N $ и $ 1 \le p  \le \infty $ через $ l_p^n $
обозначается пространство $ \R^n $ с фиксированной в нём нормой
$$
\|x\|_{l_p^n} = \begin{cases} (\sum_{j =1}^n |x_j|^p)^{1/p} \text{ при } p < \infty; \\
\max_{j =1}^n |x_j| \text{ при } p = \infty, \end{cases} x \in \R^n.
$$

Напомним ещё, что при  $ 1 \le p,q \le \infty, n \in \N $ для $ x \in \R^n $
справедливо неравенство
\begin{equation*} \tag{1.1.1}
\|x\|_{l_q^n} \le n^{(1/q -1/p)_+ } \|x\|_{l_p^n}.
\end{equation*}

Для $ d \in \N, l \in \Z_+^d $ через $ \mathcal P^{d,l} $ будем
обозначать пространство вещественных  полиномов, состоящее из всех
функций $ f: \R^d \mapsto \R $ вида
$$
f(x) = \sum_{\lambda \in \Z_+^d(l)} a_{\lambda} \cdot x^{\lambda},
x \in \R^d,
$$
где $ a_{\lambda} \in \R, x^{\lambda} = x_1^{\lambda_1} \ldots x_d^{\lambda_d},
\lambda \in \Z_+^d(l). $

При $ d \in \N, l \in \Z_+^d $ для области $ D \subset \R^d $
через $\mathcal P^{d,l}(D) $ обозначим пространство функций
$ f, $ определенных в $ D, $ для каждой из которых существует
полином $ g \in \mathcal P^{d,l} $
такой, что  сужение $ g \mid_D = f.$

Для измеримого по Лебегу множества $ D \subset \R^d $ при $ 1 \le p \le \infty $
через $ L_p(D),$ как обычно, обозначается
пространство всех вещественных измеримых на $ D $ функций $f,$
для которых определена норма
$$
\|f\|_{L_p(D)} = \begin{cases} (\int_D |f(x)|^p dx)^{1/p}, 1 \le p < \infty;\\
\supvrai_{x \in D}|f(x)|, p = \infty. \end{cases}
$$

При $ d \in \N $ для $ \lambda \in \Z_+^d $ через $ \D^\lambda $
будем обозначать оператор дифференцирования $ \D^\lambda =
\frac{\D^{|\lambda|}} {\D x_1^{\lambda_1} \ldots \D x_d^{\lambda_d}}, $
где $ |\lambda| = \sum_{j =1}^d \lambda_j. $

При обозначении известных пространств дифференцируемых функций будем
ориентироваться на [5].
Для области $ D \subset \R^d $ при $ 1 \le p \le \infty, l \in \Z_+^d $
через $ W_p^l(D) $ обозначается пространство всех функций $ f \in
L_p(D), $ для которых для каждого $ j =1,\ldots,d $  обобщенная
производная $ \D_j^{l_j} f = \frac{\D^{l_j}f }{ \D x_j^{l_j} } $
принадлежит $ L_p(D), $ с нормой
$$
\|f\|_{W_p^l(D)} = \max(\|f\|_{L_p(D)}, \max_{j =1,\ldots,d}
\|\D_j^{l_j}f\|_{L_p(D)}).
$$

Для $ x,y \in \R^d $ положим $ xy = x \cdot y = (x_1 y_1, \ldots,x_d y_d), $
а для $ x \in \R^d $ и $ A \subset \R^d $ определим
$$
x A = x \cdot A = \{xy: y \in A\}.
$$

Для $ x \in \R^d: x_j \ne 0, $ при $ j=1,\ldots,d,$ положим
$ x^{-1} = (x_1^{-1},\ldots,x_d^{-1}). $

При $ d \in \N $ для $ x \in \R^d $ положим
$$
x_+ = ((x_1)_+, \ldots, (x_d)_+),
$$
где $ t_+ = \frac{1} {2} (t +|t|), t \in \R. $

Обозначим через $ \R_+^d $ множество $ x \in \R^d, $ для которых
$ x_j >0 $ при $ j =1,\ldots,d,$ и для $ a \in \R_+^d, x \in \R^d $
положим $ a^x = a_1^{x_1} \ldots a_d^{x_d}.$

Будем также обозначать через $ \chi_A $ характеристическую функцию
множества $ A \subset \R^d. $

При $ d \in \N $ определим множества
$$
I^d = \{x \in \R^d: 0 < x_j < 1, j =1,\ldots,d\},
$$
$$
\overline I^d = \{x \in \R^d: 0 \le x_j \le 1, j =1,\ldots,d\},
$$
$$
B^d = \{x \in \R^d: -1 \le x_j \le 1, j =1,\ldots,d\}.
$$

Через $ \e $ будем обозначать вектор в $ \R^d, $ задаваемый
равенством $ \e = (1,\ldots,1). $

Напомним, что для области $ D \subset \R^d $ и вектора $ h \in \R^d $
через $ D_h $ обозначается множество
$$
D_h = \{x \in D: x +th \in D \ \forall t \in \overline I\}.
$$

Для функции $ f, $ заданной в области $ D \subset \R^d, $ и
вектора $ h \in \R^d $ определим в $ D_h $ её разность $ \Delta_h f $ с шагом
$ h, $ полагая
$$
(\Delta_h f)(x) = f(x+h) -f(x), x \in D_h,
$$
а для $ l \in \N: l \ge 2, $ в $ D_{lh} $ определим $l$-ую
разность $ \Delta_h^l f $ функции $ f $ с шагом $ h $ равенством
$$
(\Delta_h^l f)(x) = (\Delta_h (\Delta_h^{l-1} f))(x), x \in
D_{lh},
$$
положим также $ \Delta_h^0 f = f. $

Как известно, справедливо равенство
$$
(\Delta_h^l f)(\cdot) = \sum_{k=0}^l C_l^k (-1)^{l-k} f(\cdot +kh),
C_l^k = \frac{l!} {k! (l-k)!}.
$$

При $ d \in \N $ для $ j =1,\ldots,d$ через $ e_j $ будем
обозначать вектор $ e_j = (0,\ldots,0,1_j,0,\ldots,0).$

Определим пространства и классы функций, изучаемые в настоящей работе (ср. с
[5]). Но прежде введём некоторые обозначения.

Пусть $ D$ --- область в $ \R^d$ и $ 1 \le p < \infty.$
Для $ f \in L_p(D) $ при $ j =1,\ldots,d,\ l \in \Z_+ $ обозначим модуль
непрерывности в $ L_p(D) $ порядка $ l $ по $ j $-му координатному направлению через
$$
\Omega_j^l(f,t)_{L_p(D)} =
\supvrai_{h \in t B^1} \| \Delta_{h e_j}^l f\|_{L_p(D_{l h e_j})}, t \in \R_+,
$$

а также введем в рассмотрение "усредненный" модуль непрерывности
в $ L_p(D) $ порядка $ l $ по $ j $-му координатному направлению, полагая
\begin{multline*}
\Omega_j^{\prime l}(f,t)_{L_p(D)} =
\biggl((2 t)^{-1} \int_{ t B^1} \| \Delta_{\xi e_j}^l f\|_{L_p(D_{l \xi e_j})}^p d\xi\biggr)^{1 /p} = \\
\biggl((2 t)^{-1} \int_{ t B^1} \int_{D_{l \xi e_j}} |\Delta_{\xi e_j}^l f(x)|^p dx
d\xi\biggr)^{1 /p}, t \in \R_+.
\end{multline*}

Пусть теперь $ d \in \N, \alpha \in \R_+^d, 1 \le p < \infty $ и $ D $ --
область в $ \R^d. $ Тогда зададим вектор $ l = l(\alpha) \in \N^d, $ полагая
$ (l(\alpha))_j = \min \{m \in \N: \alpha_j < m \}, j =1,\ldots,d, $ а также выберем
вектор $ \l \in \Z_+^d $ такой, что $ \l < \alpha, $ и через
$ (H_p^{\alpha})^{\l}(D) $ обозначим пространство всех функций
$ f \in W_p^{\l}(D), $ для которых при $ j =1,\ldots,d $ конечна величина
$$
\sup_{t \in \R_+} t^{-(\alpha_j -\l_j)} \cdot \Omega_j^{l_j -\l_j}(\D_j^{\l_j} f, t)_{L_p(D)} < \infty,
$$
а через $ (\mathcal H_p^{\alpha})^{\l}(D) $ -- множество функций $ f \in
W_p^{\l}(D), $ для которых при $ j =1,\ldots,d $ выполняется неравенство
$$
\sup_{t \in \R_+} t^{-(\alpha_j -\l_j)} \Omega_j^{l_j -\l_j}(\D_j^{\l_j} f,
t)_{L_p(D)} \le 1.
$$

В пространстве $ (H_p^\alpha)^{\l}(D) $ вводится норма
$$
\| f \|_{(H_p^\alpha)^{\l}(D)} = \max(\| f \|_{W_p^{\l}(D)}, \max_{j =1, \ldots, d}
\sup_{t \in \R_+} t^{-(\alpha_j -\l_j)} \Omega_j^{l_j -\l_j}(\D_j^{\l_j} f, t)_{L_p(D)}).
$$

При тех же условиях на $ \alpha, p, D $ обозначим через
$ (H_p^\alpha)^\prime(D) $ пространство всех функций
$ f \in L_p(D), $ обладающих тем свойством, что для любого
$ j =1,\ldots,d $ соблюдается условие
$$
\sup_{t \in \R_+} t^{-\alpha_j} \cdot \Omega_j^{\prime l_j}(f, t)_{L_p(D)} < \infty,
$$
а через $ (\mathcal H_p^\alpha)^\prime(D) $ -- множество функций
$ f \in L_p(D), $ обладающих тем свойством, что для любого
$ j =1,\ldots,d $ соблюдается неравенство
$$
\sup_{t \in \R_+} t^{-\alpha_j} \cdot \Omega_j^{\prime l_j}(f, t)_{L_p(D)} \le 1,
\text{где} \ l = l(\alpha).
$$
В пространстве $ (H_p^\alpha)^\prime(D) $ задается норма
$$
\| f \|_{(H_p^\alpha)^\prime(D)} = \max(\| f \|_{L_p(D)}, \max_{j =1, \ldots, d}
\sup_{t \in \R_+} t^{-\alpha_j} \Omega_j^{\prime l_j}(f,t)_{L_p(D)}).
$$

Для области $ D \subset \R^d $ при $ \alpha \in \R_+^d,\ 1 \le p < \infty,
\theta \in \R: 1 \le \theta < \infty, $
полагая, как и выше, $ l = l(\alpha) \in \N^d, $ и выбирая $ \l \in \Z_+^d:
\l < \alpha, $ через $ (B_{p,\theta}^\alpha)^{\l}(D) $ обозначим пространство
всех функций $ f \in W_p^{\l}(D), $ которые для любого $ j =1,\ldots,d $
удовлетворяют условию
$$
\int_0^\infty t^{-1 -\theta (\alpha_j -\l_j)}
(\Omega_j^{l_j -\l_j}(\D_j^{\l_j}f, t)_{L_p(D)})^{\theta}\,dt < \infty,
$$
а через $ (\mathcal B_{p,\theta}^\alpha)^{\l}(D) $ обозначим множество
всех функций $ f \in W_p^{\l}(D), $ для которых при любом $ j =1,\ldots,d $
соблюдается неравенство
$$
\int_0^\infty t^{-1 -\theta (\alpha_j -\l_j)}
(\Omega_j^{l_j -\l_j}(\D_j^{\l_j}f, t)_{L_p(D)})^{\theta}\,dt \le 1.
$$

В пространстве $ (B_{p,\theta}^\alpha)^{\l}(D) $ фиксируется норма
$$
\| f \|_{(B_{p,\theta}^\alpha)^{\l}(D)} = \max\biggl(\| f \|_{W_p^{\l}(D)},
\max_{j =1,\ldots,d} \left(\int_0^\infty
t^{-1 -\theta (\alpha_j -\l_j)} (\Omega_j^{l_j -\l_j}(\D_j^{\l_j}f,
t)_{L_p(D)})^{\theta}\,dt \right)^{1/\theta}\biggr).
$$

При $ \theta = \infty $ положим $ (B_{p,\infty}^\alpha)^{\l}(D) =
(H_p^\alpha)^{\l}(D). $

При тех же условиях на параметры  обозначим через
$ (B_{p,\theta}^\alpha)^\prime(D) $ пространство всех функций
$ f \in L_p(D), $ которые при $ l = l(\alpha) $ для каждого $ j = 1,\ldots,d $
подчинены условию
$$
\int_0^\infty t^{-1 -\theta \alpha_j}
(\Omega_j^{\prime l_j}(f, t)_{L_p(D)})^{\theta}\,dt < \infty,
$$
а через $ (\mathcal B_{p,\theta}^\alpha)^\prime(D) $ -- множество функций
$ f \in L_p(D), $ которые при $ j = 1,\ldots,d $ удовлетворяют неравенству
$$
\int_0^\infty t^{-1 -\theta \alpha_j}
(\Omega_j^{\prime l_j}(f, t)_{L_p(D)})^{\theta}\,dt \le 1.
$$

В пространстве $ (B_{p,\theta}^\alpha)^\prime(D) $ определяется норма
$$
\| f \|_{(B_{p,\theta}^\alpha)^\prime(D)} = \max\biggl(\| f \|_{L_p(D)},
\max_{j =1,\ldots,d} \left(\int_0^\infty
t^{-1 -\theta \alpha_j} (\Omega_j^{\prime l_j}(f, t)_{L_p(D)})^{\theta}\,dt
\right)^{1/\theta}\biggr).
$$
При $ \theta = \infty $ положим $ (B_{p,\infty}^\alpha)^\prime(D) =
(H_p^\alpha)^\prime(D). $

В случае, когда вектор $ \l = \l(\alpha) \in \Z_+^d $ имеет компоненты
$ (\l(\alpha))_j = \max\{m \in \Z_+: m < \alpha_j\}, j =1,\ldots,d,$
пространство $ (B_{p,\theta}^\alpha)^{\l}(D) ((H_p^\alpha)^{\l}(D))$ обычно
обозначается $ B_{p,\theta}^\alpha(D) (H_p^\alpha(D)).$

Как отмечалось в [6], для $ \alpha \in \R_+^d, 1 \le p,\theta <
\infty, $ области $ D $ в $ \R^d $ имеют место соотношения
\begin{multline*} \tag{1.1.2}
(B_{p, \theta}^\alpha)^\prime(D) \subset
(H_p^\alpha)^\prime(D);
(\mathcal B_{p, \theta}^\alpha)^\prime(D) \subset
c_1(\alpha) (\mathcal H_p^\alpha)^\prime(D); \\
\| f\|_{(H_p^\alpha)^\prime(D)} \le c_1(\alpha)
\| f\|_{(B_{p, \theta}^\alpha)^\prime(D)},
\end{multline*}
где $ c_1(\alpha) = \max_{j =1,\ldots,d} 2^{2+\alpha_j}, $
а также
\begin{equation*} \tag{1.1.3}
(B_{p, \theta}^\alpha)^{\l}(D) \subset (B_{p, \theta}^\alpha)^\prime(D), \\
(\mathcal B_{p, \theta}^\alpha)^{\l}(D) \subset (\mathcal B_{p, \theta}^\alpha)^\prime(D)
\end{equation*}
и
\begin{multline*} \tag{1.1.4}
\| f\|_{(B_{p, \theta}^\alpha)^\prime(D)} \le
\| f\|_{(B_{p, \theta}^\alpha)^{\l}(D)},
\alpha \in \R_+^d, \l \in \Z_+^d: \l < \alpha, 1 \le p < \infty,
1 \le \theta \le \infty, \\
D \text{ -- произвольная область в } \R^d.
\end{multline*}

Обозначим через $ C^\infty(D) $ пространство бесконечно дифференцируемых
функций в области $ D \subset \R^d, $ а через $ C_0^\infty(D) $ -- пространство
функций $ f \in C^\infty(\R^d), $ у которых носитель $ \supp f \subset D. $

В заключение этого пункта введём ещё несколько обозначений.
Для линейного нормированного пространства $ X $ (над $ \R$) обозначим
$ B(X) = \{x \in X: \|x\|_X \le 1\}. $

Для линейных нормированных пространств $ X,Y $ через $ \mathcal B(X,Y) $
обозначим линейное нормированное пространство, состоящее из непрерывных
линейных операторов $ T: X \mapsto Y, $ с нормой
$$
\|T\|_{\mathcal B(X,Y)} = \sup_{x \in B(X)} \|Tx\|_Y.
$$
\bigskip

1.2. В этом пункте будут построены семейства операторов проектирования
на подпространства кусочно-полиномиальных функций, и описаны их
свойства, которые используются при доказательстве основных результатов
работы.

Для $ d \in \N, x, y \in \R^d $ будем писать $ x \le y ( x < y), $
если для каждого $ j =1,\ldots,d $ выполняется неравенство $ x_j
\le y_j (x_j < y_j). $

Для $ d \in \N, m, n \in \Z^d: m \le n, $ обозначим
$$
\Nu_{m, n}^d = \{ \nu \in \Z^d: m \le \nu \le n \} = \prod_{j=1}^d
\Nu_{m_j,n_j}^1.
$$

При $ d \in \N $ для $ t \in \R^d $ через $ 2^{t} $ будем обозначать
вектор $ 2^{t} = (2^{t_1}, \ldots, 2^{t_d}). $

Для $ d \in \N, \kappa \in \Z^d, \nu \in \Z^d $ обозначим через
$ \chi_{\kappa, \nu}^d $ характеристическую функцию множества
$ Q_{\kappa, \nu}^d = 2^{-\kappa} \nu +2^{-\kappa} I^d. $ Понятно,
что при $ d \in \N, \kappa \in \Z^d, \nu \in \Z^d $ имеют место
равенства
$$
Q_{\kappa, \nu}^d = \prod_{j=1}^d
Q_{\kappa_j,\nu_j}^1, \\
\chi_{\kappa, \nu}^d(x) = \prod_{j=1}^d
\chi_{\kappa_j,\nu_j}^1(x_j), x \in \R^d.
$$

Введем в рассмотрение пространства кусочно-полиномиальных функций.

Для $ d \in \N, l \in \Z_+^d $ и $ \kappa \in \Z_+^d$ через
$ \mathcal P_{\kappa}^{d, l}$ обозначим линейное подпространство в
$ L_\infty(I^d),$ состоящее из функций $ f \in L_\infty(I^d), $ для
каждой из которых существует набор  полиномов $ \{f_{\nu} \in
\mathcal P^{d, l}, \nu \in \Nu_{0, 2^{\kappa} -\e}^d\} $ такой, что
\begin{equation*} \tag{1.2.1}
f = \sum_{ \nu \in \Nu_{0, 2^{\kappa} -\e}^d} f_{\nu} \chi_{\kappa, \nu}^d.
\end{equation*}

Следующая лемма взята из [7].

Лемма 1.2.1

Пусть $ d \in \N $ и $ \kappa, \kappa^\prime, \nu, \nu^\prime \in \Z^d $
таковы, что $ \kappa^\prime \le \kappa, $ а $ Q_{\kappa, \nu}^d \cap
Q_{\kappa^\prime, \nu^\prime}^d \ne \emptyset. $
Тогда имеет место включение
\begin{equation*} \tag{1.2.2}
Q_{\kappa, \nu}^d \subset Q_{\kappa^\prime, \nu^\prime}^d.
\end{equation*}

Из (1.2.1), (1.2.2) следует, что при $ d \in \N, l \in \Z_+^d $ для $ \kappa,
\kappa^\prime \in \Z_+^d: \kappa^\prime \le \kappa, $ справедливо
включение
\begin{equation*} \tag{1.2.3}
\mathcal P_{\kappa^\prime}^{d, l} \subset \mathcal P_{\kappa}^{d, l}.
\end{equation*}

При $ k \in \Z_+, \alpha \in \R_+^d $ определим $ \kappa(k,\alpha) $
как вектор, имеющий компоненты
\begin{equation*}
(\kappa(k,\alpha))_j = [k /\alpha_j], j =1,\ldots,d, \ ([a] \text{ -- целая часть} \ a),
\end{equation*}
и для $ d \in \N, l \in \Z_+^d, \alpha \in \R_+^d $ и $ k \in \Z_+ $
через $ \mathcal P_k^{d,l,\alpha} $ обозначим линейное подпространство в
$ L_\infty(I^d)$, определяемое равенством
$$
\mathcal P_k^{d,l,\alpha} = \mathcal P_{\kappa}^{d, l}
$$
при $ \kappa = \kappa(k,\alpha). $

Из (1.2.3) вытекает, что для $ d \in \N, l \in \Z_+^d, \alpha \in \R_+^d, k \in \N $
справедливо включение
\begin{equation*} \tag{1.2.4}
\mathcal P_{k -1}^{d,l,\alpha} \subset \mathcal P_k^{d,l,\alpha}.
\end{equation*}

Для $ x,y \in \R^d $ будем обозначать $ (x,y) = \sum_{j=1}^d x_j y_j. $

Обозначая через $ R_k^{d,l,\alpha} = \dim \mathcal P_k^{d,l,\alpha} $,
отметим (см. [8]), что для $ d \in \N, l \in \Z_+^d, \alpha \in \R_+^d $
существуют константы $ c_1(d,l,\alpha) >0 $ и $ c_2(d,l,\alpha)>0 $ такие, что
при $ k \in \Z_+ $ выполняется неравенство
\begin{equation*} \tag{1.2.5}
c_1 2^{k(\alpha^{-1},\e)} < R_k^{d,l,\alpha} < c_2 2^{k(\alpha^{-1},\e)}.
\end{equation*}
Для $ d \in \N, l \in \Z_+^d, \alpha \in\R_+^d $ и $ k \in \Z_+ $  обозначим
через $ \I_k^{d,l,\alpha}: \mathcal P_k^{d,l,\alpha} \mapsto \R^{R_k^{d,l,\alpha}} $
линейный изоморфизм, определяемый для $ f \in \mathcal P_k^{d,l,\alpha}  $

равенством
$$
\I_k^{d,l,\alpha} f = \{ f_\nu( 2^{-\kappa} \nu +2^{-\kappa} \lambda): \nu
\in \Nu_{0,2^\kappa -\e}^d, \lambda \in \Z_+^d(l) \},
$$
где $ \kappa = \kappa(k,\alpha), $ а $ \{f_\nu \in \mathcal P^{d,l}, \nu \in
\Nu_{0, 2^\kappa -\e}^d\} $ -- набор
полиномов, удовлетворяющих (1.2.1).

Из [8] можно извлечь лемму.

Лемма 1.2.2

Пусть $ d \in \N, l \in \Z_+^d, \alpha \in \R_+^d, 1 \le p \le \infty. $
Тогда при любом $ k \in \Z_+ $ линейный изоморфизм $ \I_{k}^{d, l,\alpha} $
подпространства $ \mathcal P_{k}^{d, l,\alpha} $ на пространство
$ \R^{R_{k}^{d, l,\alpha}}, $ обладает тем свойством, что
для $ f \in \mathcal P_k^{d,l,\alpha} $ соблюдаЮтся неравенства
\begin{equation*} \tag{1.2.6}
c_3 \|f\|_{L_p(I^d)} \le 2^{-k(\alpha^{-1}, \e) p^{-1}}
\| \I_k^{d,l,\alpha} f \|_{l_p^{R_k^{d,l,\alpha}}} \le c_4
\|f\|_{L_p(I^d)}
\end{equation*}
с некоторыми константами $ c_3 >0, c_4 >0, $ зависящими только
от $ d, l, \alpha, p. $

При определении операторов проектирования на подпространства
$ \mathcal P_{\kappa}^{d, l} $ используются операторы из следующего
предложения.

Предложение 1.2.3

Пусть $ d \in \N, l \in \Z_+^d. $ Тогда
для любого $ \delta \in \R_+^d $ и $ x^0 \in \R^d $ для $ Q = x^0 +\delta I^d $
существует единственный линейный оператор
$ P_{\delta, x^0}^{d,l}: L_1(Q) \mapsto \mathcal P^{d, l}, $
обладающий следующими свойствами:

1) для $ f \in \mathcal P^{d, l} $ имеет место равенство
\begin{equation*}
P_{\delta, x^0}^{d, l}(f \mid_Q) = f,
\end{equation*}

2)
\begin{equation*}
\Ker P_{\delta, x^0}^{d, l} = \biggl\{\,f \in L_1(Q):
\int\limits_{Q} f(x) g(x) \,dx =0\ \forall g \in \mathcal P^{d, l}\,\biggr\},
\end{equation*}

причем, существуют константы $ c_5(d, l) >0 $ и $ c_6(d, l) >0 $
такие, что

3) при $ 1 \le p \le \infty $ для $ f \in L_p(Q) $ справедливо неравенство
\begin{equation*}
\|P_{\delta, x^0}^{d, l} f \|_{L_p(Q)} \le c_5 \|f\|_{L_p(Q)},
\end{equation*}

4) при $ 1 \le p < \infty $ для $ f \in L_p(Q) $ выполняется неравенство
\begin{equation*}
\| f -P_{\delta, x^0}^{d, l} f \|_{L_p(Q)} \le c_6 \sum_{j =1}^d
\delta_j^{-1/p} \biggl(\int_{\delta_j B^1} \int_{Q_{(l_j +1) \xi e_j}}
|\Delta_{\xi e_j}^{l_j +1} f(x)|^p dx d\xi\biggr)^{1/p}.
\end{equation*} (см. [7])

Для $ d \in \N, l, \kappa \in \Z_+^d, \nu \in \Nu_{0, 2^{\kappa} -\e}^d $
определим непрерывный линейный оператор $ S_{ \kappa, \nu}^{d, l}: L_1(I^d)
\mapsto \mathcal P^{d, l}(I^d) \cap L_\infty(I^d), $ полагая для $ f \in
L_1(I^d) $ значение
$$
S_{\kappa, \nu}^{d, l} f =
(P_{\delta, x^0}^{d, l} (f \mid_{(x^0 +\delta I^d)})) \mid_{I^d}
$$
при $ \delta = 2^{-\kappa}, x^0 = 2^{-\kappa} \nu $ (см. предложение 1.2.3).

Определим при $ d \in \N, l, \kappa \in \Z_+^d $ линейный непрерывный оператор
$ E_{\kappa}^{d, l}: L_1(I^d) \mapsto \mathcal P_{\kappa}^{d, l} \cap L_\infty (I^d) $
равенством
$$
E_{\kappa}^{d, l} f = \sum_{\nu \in \Nu_{0, 2^{\kappa} -\e}^d}
(S_{\kappa, \nu}^{d, l} f) \chi_{\kappa, \nu}^d, f \in L_1(I^d).
$$

Отметим некоторые свойства этих операторов, установленные в [7].

Лемма 1.2.4

Пусть  $ d \in \N, l \in \Z_+^d.$ Тогда
справедливы следующие утверждения:

1) при $ \kappa \in \Z_+^d $ для $ f \in \mathcal P_{\kappa}^{d, l} $
соблюдается равенство
\begin{equation*}
E_{\kappa}^{d, l} f = f;
\end{equation*}

2) при $ \kappa \in \Z_+^d $ ядро
\begin{equation*}
\Ker E_{\kappa}^{d, l} = \biggl\{ \, f \in L_1(I^d):
\int\limits_{I^d} f(x) g(x) \,dx =0 \ \forall g \in
\mathcal P_{\kappa}^{d, l} \,\biggr\};
\end{equation*}

3) существует константа $ c_7(d,l) >0 $ такая, что при $ \kappa \in \Z_+^d $
и $ 1 \le q \le \infty $ для $ f \in L_q(I^d) $ выполняется неравенство
\begin{equation*}
\| E_{\kappa}^{d, l} f \|_{L_q(I^d)} \le c_7 \| f \|_{L_q(I^d)}.
\end{equation*}

Для $ d  \in \N, l \in \Z_+^d, \alpha \in \R_+^d $ при $ k \in \Z_+ $ положим
$$
E_k^{d,l,\alpha} = E_{\kappa(k,\alpha)}^{d,l}, \\
E_{-1}^{d,l,\alpha} =0.
$$

Частным случаем леммы 1.2.4 является следующая лемма.

Лемма 1.2.5

Пусть  $ d \in \N, l \in \Z_+^d, \alpha \in \R_+^d. $ Тогда имеют место
следующие утверждения:

1) при $ k \in \Z_+ $ для $ f \in \mathcal P_k^{d, l,\alpha} $
соблюдается равенство
\begin{equation*} \tag{1.2.7}
E_k^{d, l,\alpha} f = f;
\end{equation*}

2) при $ k \in \Z_+ $ ядро
\begin{equation*} \tag{1.2.8}
\Ker E_k^{d, l,\alpha} = \biggl\{ \, f \in L_1(I^d):
\int\limits_{I^d} f(x) g(x) \,dx =0 \ \forall g \in
\mathcal P_k^{d, l,\alpha} \,\biggr\};
\end{equation*}

3) существует константа $ c_7(d,l) >0 $ такая, что при $ k \in \Z_+ $
и $ 1 \le q \le \infty $ для $ f \in L_q(I^d) $ выполняется неравенство
\begin{equation*} \tag{1.2.9}
\| E_k^{d, l,\alpha} f \|_{L_q(I^d)} \le c_7 \| f \|_{L_q(I^d)}.
\end{equation*}

Лемма 1.2.6

Пусть $ d \in \N, \alpha \in \R_+^d, 1 \le p < \infty, 1 \le q \le \infty $
удовлетворяют условию
\begin{equation*} \tag{1.2.10}
1 -(\alpha^{-1},\e)(p^{-1} -q^{-1})_+  >0
\end{equation*}
и $ l = l(\alpha). $
Тогда существует константа $ c_8(d,\alpha,p,q) >0 $ такая, что для  любой
функции $ f \in B((H_p^\alpha)^\prime(I^d)) $ при $ k \in \Z_+ $
соблюдается неравенство
\begin{equation*} \tag{1.2.11}
\| f -E_k^{d,l -\e,\alpha} f \|_{L_q(I^d)} \le
c_8 2^{-k(1 -(\alpha^{-1},\e)(p^{-1} -q^{-1})_+ )}.
\end{equation*}

Доказательство.

В условиях леммы для $ f \in B((H_p^\alpha)^\prime(I^d)) $ при $ k \in \Z_+,
l = l(\alpha) \in \N^d $ в силу неравенства (2.1.24) из [6] и
леммы 4.1.2 из [8] существует $ g \in \mathcal P_k^{d, l -\e,\alpha}, $ для
которого выполняется неравенство
\begin{equation*}
\| f -g \|_{L_q(I^d)} \le c_9 2^{-k(1 -(\alpha^{-1},\e)(p^{-1} -q^{-1})_+ )}
\end{equation*}
с константой $ c_9(d,\alpha,p,q) >0, $
а, следовательно, благодаря (1.2.7), (1.2.9), имеет место оценка
\begin{multline*}
\| f -E_k^{d,l -\e,\alpha} f \|_{L_q(I^d)} =
\| f -g +g -E_k^{d,l -\e,\alpha} f \|_{L_q(I^d)} = \\
\| f -g +E_k^{d,l -\e,\alpha} g -E_k^{d,l -\e,\alpha} f \|_{L_q(I^d)} =
\| f -g +E_k^{d,l -\e,\alpha}(g -f) \|_{L_q(I^d)} \le \\
\| f -g \|_{L_q(I^d)} +\| E_k^{d,l -\e,\alpha}(g -f) \|_{L_q(I^d)} \le
(1 +c_7) \| f -g \|_{L_q(I^d)} \le \\
c_8 2^{-k(1 -(\alpha^{-1},\e)(p^{-1} -q^{-1})_+ )}. \square
\end{multline*}

Лемма 1.2.7

При $ d \in \N, l \in \Z_+^d, \alpha \in \R_+^d $ для $ j,k \in \Z_+: j \le k, $
справедливы равенства
\begin{equation*} \tag{1.2.12}
E_j^{d,l,\alpha} E_k^{d,l,\alpha} = E_k^{d,l,\alpha} E_j^{d,l,\alpha} =
E_j^{d,l,\alpha}.
\end{equation*}

Доказательство.

Поскольку ввиду (1.2.4) при $ j \le k$ для $ f \in L_1(I^d)$ имеет место
включение
$ E_j^{d,l,\alpha} f \in \mathcal P_j^{d,l,\alpha} \subset
\mathcal P_k^{d,l,\alpha} $, то из (1.2.7) следует, что
$ E_k^{d,l,\alpha} E_j^{d,l,\alpha} f = E_j^{d,l,\alpha} f.$

Далее, при $ n \in \Z_+$ из (1.2.7) вытекает, что $ E_n^{d,l,\alpha} $ суть
оператор проектирования пространства $ L_1(I^d) $ на подпространство
$ \mathcal P_n^{d,l,\alpha}$ параллельно подпространству $ \Ker E_n^{d,l,\alpha}$.
Поэтому, учитывая, что в силу (1.2.8), (1.2.4) при $ j \le k $ ядро
$ \Ker E_k^{d,l,\alpha} \subset \Ker E_j^{d,l,\alpha} $, для $ f \in
L_1(I^d) $ имеет место представление
\begin{equation*}
f = E_k^{d,l,\alpha} f +g,\ g \in \Ker E_k^{d,l,\alpha} \subset \Ker E_j^{d,l,\alpha} ,\\
E_k^{d,l,\alpha} f = E_j^{d,l,\alpha} E_k^{d,l,\alpha} f +h,\ h \in \Ker E_j^{d,l,\alpha} ,
\end{equation*}
т.е.
$$
f = E_j^{d,l,\alpha} E_k^{d,l,\alpha} f +(h +g),
$$
где $ (h +g) \in \Ker E_j^{d,l,\alpha} ,\ E_j^{d,l,\alpha} E_k^{d,l,\alpha} f
\in \mathcal P_j^{d,l,\alpha} $.
А это значит, что $ E_j^{d,l,\alpha} E_k^{d,l,\alpha} f = E_j^{d,l,\alpha} f. \square$

Лемма 1.2.8

При $ d \in \N, l \in \Z_+^d, \alpha \in \R_+^d $ для операторов
$ \mathcal E_k^{d,l,\alpha} = E_k^{d,l,\alpha} -E_{k -1}^{d,l,\alpha},\ k \in \Z_+,\
E_{-1}^{d,l,\alpha} =0, $ выполняются равенства
\begin{multline*} \tag{1.2.13}
\mathcal E_j^{d,l,\alpha} \mathcal E_k^{d,l,\alpha} = \begin{cases}
\mathcal E_k^{d,l,\alpha} ,\text{ при } j = k; \\
0,\text{ при } j \ne k.
\end{cases}
\end{multline*}

Доказательство.

Пусть $ j = k. $ Тогда в силу (1.2.12)
\begin{multline*}
(\mathcal E_k^{d,l,\alpha})^2 = (E_k^{d,l,\alpha})^2 -E_{k -1}^{d,l,\alpha}
E_k^{d,l,\alpha} -E_k^{d,l,\alpha} E_{k -1}^{d,l,\alpha} +(E_{k -1}^{d,l,\alpha})^2 \\
= E_k^{d,l,\alpha} -E_{k -1}^{d,l,\alpha} -E_{k -1}^{d,l,\alpha} +
E_{k -1}^{d,l,\alpha} = \mathcal E_k^{d,l,\alpha}.
\end{multline*}

Пусть $ j \ne k. $ Предположим, что $ j < k. $ Тогда, снова используя (1.2.12),
выводим
\begin{multline*}
\mathcal E_j^{d,l,\alpha} \mathcal E_k^{d,l,\alpha} = E_j^{d,l,\alpha}
E_k^{d,l,\alpha} -E_{j -1}^{d,l,\alpha} E_k^{d,l,\alpha} -E_j^{d,l,\alpha}
E_{k -1}^{d,l,\alpha} +E_{j -1}^{d,l,\alpha} E_{k -1}^{d,l,\alpha}  \\
= E_j^{d,l,\alpha} -E_{j -1}^{d,l,\alpha} -E_j^{d,l,\alpha} +E_{j -1}^{d,l,\alpha} =0.
\end{multline*}

Аналогично проверяется (1.2.13) при $ j > k. \square$

Лемма 1.2.9

При $ d \in \N, l \in \Z_+^d, \alpha \in \R_+^d $ для $ k \in \Z_+,\ j \in \N$
определим оператор $ \mathfrak E_{k,k +j}^{d,l,\alpha} = E_{k +j}^{d,l,\alpha}
-E_k^{d,l,\alpha} $, а через $ \mathfrak P_{k,k +j}^{\prime\,d,l,\alpha} $
обозначим его образ
$ \mathfrak P_{k,k +j}^{\prime\,d,l,\alpha} = \operatorname{Im}
\mathfrak E_{k,k +j}^{d,l,\alpha} $.

Тогда для $ k \in \Z_+,\ j \in \N $ справедливы следующие соотношения:

1)
\begin{equation*} \tag{1.2.14}
\mathfrak E_{k,k +j}^{d,l,\alpha} f = f \text{ для } f \in \mathfrak P_{k,k +j}^{\prime\,d,l,\alpha};
\end{equation*}

2)
\begin{gather*}  \tag{1.2.15}
\mathcal P_{k +j}^{d,l,\alpha} = \mathcal P_k^{d,l,\alpha} +
\mathfrak P_{k,k +j}^{\prime\,d,l,\alpha}, \\
\text{ причём } \\
\mathfrak P_{k,k +j}^{\prime\,d,l,\alpha} =
\{f \in \mathcal P_{k +j}^{d,l,\alpha}: E_k^{d,l,\alpha} f = 0\} \\
\text{ и } \\
\mathcal P_k^{d,l,\alpha} \bigcap \mathfrak P_{k,k +j}^{\prime\,d,l,\alpha} = \{0\}.
\end{gather*}

Доказательство.

Сначала проверим (1.2.14). Пусть $ f \in\mathfrak P_{k,k +j}^{\prime\,d,l,\alpha} $.
Тогда $ f = \mathfrak E_{k,k +j}^{d,l,\alpha} g $, где $ g \in L_1(I^d)$. При
этом в силу (1.2.12) имеем
\begin{multline*}
\mathfrak E_{k,k +j}^{d,l,\alpha} f
= (\mathfrak E_{k,k +j}^{d,l,\alpha} )^2 g = (E_{k +j}^{d,l,\alpha})^2 g -
E_{k +j}^{d,l,\alpha} E_k^{d,l,\alpha} g-E_k^{d,l,\alpha} E_{k +j}^{d,l,\alpha} g
+(E_k^{d,l,\alpha})^2 g \\
= E_{k +j}^{d,l,\alpha} g -E_k^{d,l,\alpha} g -E_k^{d,l,\alpha} g +
E_k^{d,l,\alpha} g = \mathfrak E_{k,k +j}^{d,l,\alpha} g = f.
\end{multline*}

Теперь убедимся в справедливости (1.2.15). В виду (1.2.7) для $ f \in
\mathcal P_{k +j}^{d,l,\alpha} $ имеем
\begin{equation*}
f = E_{k +j}^{d,l,\alpha} f = \mathfrak E_{k,k +j}^{d,l,\alpha} f +
E_k^{d,l,\alpha} f
\end{equation*}
и, значит,
$$
\mathcal P_{k +j}^{d,l,\alpha} = \mathcal P_k^{d,l,\alpha} +
\mathfrak P_{k,k +j}^{\prime\,d,l,\alpha}.
$$

Далее, пусть $ f \in \mathfrak P_{k,k +j}^{\prime\,d,l,\alpha}.$ Тогда
существует $ g \in L_1(I^d) $, для которого
$$
f = \mathfrak E_{k,k +j}^{d,l,\alpha} g = E_{k +j}^{d,l,\alpha} g -E_k^{d,l,\alpha} g,
$$
откуда, учитывая (1.2.4), заключаем, что $ f \in \mathcal P_{k +j}^{d,l,\alpha}, $
и, используя (1.2.12), получаем
\begin{equation*}
E_k^{d,l,\alpha} f = E_k^{d,l,\alpha} \mathfrak E_{k,k +j}^{d,l,\alpha} g =
(E_k^{d,l,\alpha} E_{k +j}^{d,l,\alpha}) g -(E_k^{d,l,\alpha})^2 g =
E_k^{d,l,\alpha} g -E_k^{d,l,\alpha} g =0.
\end{equation*}

Обратно, пусть $ f \in \mathcal P_{k +j}^{d,l,\alpha} $ и
$ E_k^{d,l,\alpha} f = 0 $. Тогда, принимая во внимание (1.2.7), имеем
$$
f = E_{k +j}^{d,l,\alpha} f = E_{k +j}^{d,l,\alpha} f -E_k^{d,l,\alpha} f = \mathfrak E_{k,k +j}^{d,l,\alpha} f \in
\mathfrak P_{k,k +j}^{\prime\,d,l,\alpha}.
$$

Наконец, если $ f \in \mathcal P_k^{d,l,\alpha} \bigcap \mathfrak P_{k,k +j}^{\prime\,d,l,\alpha} $,
то, применяя (1.2.7) с учетом сказанного выше, находим, что
$ f = E_k^{d,l,\alpha} f =0. \square$
\bigskip

1.3. В этом пункте мы напомним определение $n$-поперечника по
Бернштейну выпуклых симметричных множеств, употребляемое в [1], и
приведём сведения о них, необходимые для вывода основного
результата этой работы.

Пусть $C$ --- выпуклое симметричное (относительно $0$) подмножество
линейного нормированного пространства $X$ над $\R,\ n \in \N$ и $ \mathcal M_n(X)$ ---
совокупность всех подпространств $ L \subset X$, у которых размерность
$ \dim L = n$, а $ B(X) = \{\,x \in X: \|x\|_X \le 1\,\}$ --- единичный шар
в $X$. Тогда $n$-поперечником по Бернштейну множества $C$ в пространстве
$X$ называется величина
$$
b_n(C,X) = \sup\limits_{ L \in \mathcal M_n(X)} \sup\{\,\epsilon \ge 0:
(\epsilon \cdot B(X)) \cap L \subset C\,\}.
$$

Отметим некоторые свойства $n$-поперечника по Бернштейну, приведенные в
[2].

Предложение 1.3.1

1) Для $ a > 0$ и выпуклого симметричного множества $ C \subset X$
при $ n \in \N$ верно равенство
\begin{equation*} \tag{1.3.1}
b_n(a \cdot C,X) = a \cdot b_n(C,X).
\end{equation*}

2) Для выпуклых симметричных подмножеств $ C_1 \subset C_2$ пространства
$X$ при $ n \in \N$ соблюдается неравенство
\begin{equation*} \tag{1.3.2}
b_n(C_1, X) \le b_n(C_2, X).
\end{equation*}

3) Пусть $ X, Y$ --- линейные нормированные пространства над $\R$,
$C$ --- выпуклое симметричное подмножество $X$ и $U: X \mapsto Y$ ---
непрерывный линейный оператор, у которого $ \Ker U = \{0\}.$ Тогда при
$ n \in \N$ имеет место неравенство
\begin{equation*} \tag{1.3.3}
b_n(U(C), Y) \le \|U\|_{\mathcal B(X,Y)} b_n(C, X).
\end{equation*}

4) Для выпуклого симметричного подмножества $C$ пространства $X$
для $ n \in \N$ выполняется неравенство
\begin{equation*} \tag{1.3.4}
b_{n +1}(C,X) \le b_n(C,X).
\end{equation*}

5) Пусть $ X$ --- линейное нормированное пространство над $\R$, $ Y \subset X $ ---
линейное подпространство в $X, $ а $ Y \cap X $ --- подпространство $ Y $ с
фиксированной в нём нормой, индуцированнной из $ X. $ Пусть ещё
$C$ --- выпуклое симметричное подмножество $ Y.$ Тогда при $ n \in \N$
соблюдается равенство
\begin{equation*} \tag{1.3.5}
b_n(C, Y \cap X) = b_n(C, X).
\end{equation*}

Доказательство.

Соотношения (1.3.1) --- (1.3.4) установлены в [2]. Убедимся в
справедливости равенства (1.3.5).
В условиях п. 5) сначала пусть $ L \in \mathcal M_n(Y) $ и $ \epsilon >0 $
таковы, что имеет место включение
$$
L \cap (\epsilon B(Y \cap X)) \subset C.
$$
Тогда $ L \in \mathcal M_n(X) $ и
\begin{multline*}
C \supset L \cap (\epsilon B(Y \cap X)) = L \cap (\epsilon (Y \cap B(X))) =
L \cap (Y \cap (\epsilon B(X))) = \\
 (L \cap Y) \cap (\epsilon B(X)) = L \cap (\epsilon B(X)),
\end{multline*}
а, значит, $ \epsilon \le b_n(C,X), $ и
\begin{equation*}
b_n(C, Y \cap X) = \sup_{L \in \mathcal M_n(Y), \epsilon \ge 0: L \cap (\epsilon
B(Y \cap X)) \subset C} \epsilon \le b_n(C, X).
\end{equation*}

Теперь пусть $ L \in \mathcal M_n(X) $ и $ \epsilon >0 $ таковы, что справедливо
включение
$$
L \cap (\epsilon B(X)) \subset C.
$$
При этом для любого $ x \in L: x \ne 0, $ элемент
$$
(\epsilon / \|x\|_X) x \in L \cap (\epsilon B(X)) \subset C \subset Y,
$$
а, следовательно, $ x \in Y, $ и, значит, $ L \subset Y, $ т.е. $ L \in
\mathcal M_n(Y), $ и соблюдается включение
\begin{equation*}
L \cap (\epsilon B(Y \cap X)) = L \cap (\epsilon (Y \cap B(X))) \subset
L \cap (\epsilon B(X)) \subset C,
\end{equation*}
что влечёт неравенство $ \epsilon \le b_n(C, Y \cap X). $
Таким образом, получаем
\begin{equation*}
b_n(C, X) = \sup_{L \in \mathcal M_n(X), \epsilon \ge 0: L \cap (\epsilon
B(X)) \subset C} \epsilon \le b_n(C, Y \cap X).
\end{equation*}

Сопоставляя сказанное, приходим к (1.3.5). $\square$

Для выпуклого симметричного подмножества $C$ линейного пространства $X$
над $\R$ через $\mu_C$ обозначим его функционал Минковского, т.е.
функционал, определяемый равенством

$$
\mu_C(x) = \inf\{\,t>0: t^{-1} x \in C\,\}.
$$

Нам будут полезны следующие хорошо известные свойства функционала
Минковского.

Предложение 1.3.2

Пусть $ C $ -- выпуклое симметричное подмножество линейного пространства $X.$
Тогда

1) для $ x \in C$ имеет место неравенство
$$
\mu_C(x) \le 1;
$$

2) если для $ x \in X$ значение $ \mu_C(x) < 1$, то $ x \in C$;

3) если $ C = \bigcap \limits_{j =1}^{j_0} C_j$, где
$ C_j$ --- выпуклое симметричное множество в $X,\ j =1,\ldots,j_0$, то
\begin{equation*} \tag{1.3.6}
\mu_C(x) = \max_{j =1,\ldots,j_0} \mu_{C_j}(x).
\end{equation*}

При выводе оценки сверху поперечника $b_n(B((H_p^\alpha)^\prime(I^d)),
L_q(I^d)) $ существенную роль играет установленная автором в [2]

Теорема 1.3.3

Пусть $C$ --- выпуклое симметричное подмножество линейного нормированного
пространства $X$ над $\R$. Тогда для $ n \in \N$ справедливо равенство
\begin{equation*} \tag{1.3.7}
b_n(C, X) = \sup\limits_{L \in \mathcal M_n(X)} \inf\limits_{x \in L \setminus
\{0\}} \frac{ \|x\|_X } { \mu_C(x)}.
\end{equation*}

Приведём необходимые нам соотношения для бернштейновских поперечников
эллипсоидов, вытекающие из результатов, полученных в [9], [10], [1].
Но прежде для $ N \in \N,\ 1 \le p \le \infty$ и $\rho = (\rho_1,\ldots,\rho_N):
\rho_j \ge 0,\ j =1,\ldots,N$ через $ B_p^N(\rho)$ обозначим множество
\begin{multline*}
B_p^N(\rho) = \biggl\{\,x \in \R^N: \sum_{j =1}^N \biggl(\frac{ |x_j| }
{\rho_j} \biggr)^p \le 1\,\biggr\},\text{ при } 1 \le p < \infty,\\
B_\infty^N(\rho) = \{\,x \in \R^N: |x_j|\le \rho_j,\ j = 1,\ldots,N\,\}.
\end{multline*}

Теорема 1.3.4

Пусть $ n, N \in \N: n \le N,$ и $ \rho = (\rho_1,\ldots,\rho_N): \rho_j
\ge \rho_{j +1} \ge 0,\ j = 1,\ldots,N -1$. Тогда
\begin{equation*} \tag{1.3.8}
b_n(B_\infty^N(\rho), l_2^N) \le \biggl(2 \cdot \negthickspace
\sum_{\{\,j: \frac {n} {2} < j \le N\,\} } \frac {\rho_j^2} {n} \biggr)^{1/2}.
\end{equation*}

Теорема 1.3.5

Пусть $ n, N \in \N: n \le N,\ 1 \le p < q \le \infty$ и
$ \rho = (\rho_1,\ldots,\rho_N): \rho_j \ge \rho_{j +1} > 0,\ j =1,\ldots,N -1$.
Тогда
\begin{equation*} \tag{1.3.9}
b_n(B_p^N(\rho), l_q^N) = \biggl( \sum_{j =1}^n
\bigl(\rho_j\bigr)^{\frac{pq} {p -q}} \biggr)^{\frac{p -q} {pq}}, q < \infty,
\end{equation*}
\begin{equation*} \tag{1.3.9}^\prime
b_n(B_p^N(\rho), l_\infty^N) =
\biggl( \sum_{j =1}^n (\rho_j)^{-\!p} \biggr)^{\!-1/p}.
\end{equation*}

В заключение приведём некоторые полезные для нас факты, касающиеся
колмогоровского и бернштейновского поперечников.

Из результатов [11] следует, что при $n \in \N,\ 1 \le p,q \le \infty$
имеет место равенство
$$
b_n(B(l_p^{2n}), l_q^{2n}) = \frac{1} {d_n(B(l_{p^\prime}^{2n}), l_{q^\prime}^{2n})},
\text{ где } p^\prime = \frac {p} {p -1}, \ q^\prime = \frac {q} {q -1},
$$
что в соединении с известными оценками колмогоровских поперечников шаров (см.
[12], [13]), в частности, дает
\begin{equation*} \tag{1.3.10}
b_n(B(l_p^{2n}), l_q^{2n}) \ge c_1 \begin{cases}
n^{(1/q -1/p)}, p < q, q \le p \le 2, \\
1, 2 \le q \le p, \\
n^{(1/q -1/2)}, q \le 2 \le p.
\end{cases}
\end{equation*}
\bigskip

1.4. В этом пункте приведем используемые в следующем параграфе соотношения
между нормами образов и прообразов при некоторых отображениях рассматриваемых
пространств.

При $ d \in \N $ для $ \delta \in \R_+^d $ и $ x^0 \in \R^d $ обозначим через
$ h_{\delta, x^0} $ отображение, которое каждой функции $ f, $ заданной на
некотором множестве $ S \subset \R^d, $  ставит в соответствие функцию
$ h_{\delta, x^0} f, $ определяемую на множестве $ \{ x \in \R^d: x^0 +\delta
x \in S\} = \delta^{-1} (S -x^0) $ равенством $ (h_{\delta, x^0} f)(x) =
f(x^0 +\delta x). $ Так как для $ \delta \in \R_+^d, x^0 \in \R^d $
отображение $ \R^d \ni x \mapsto x^0 +\delta x \in \R^d $ ---
взаимно однозначно, то отображение $ h_{\delta, x^0} $ является
биекцией на себя  множества всех функций с областью определения
в $ \R^d. $ При этом обратное  отображение $ h_{\delta, x^0}^{-1} $
задается равенством
\begin{equation*}
(h_{\delta, x^0}^{-1} f)(x) = f(\delta^{-1} (x -x^0)) =
(h_{\delta^\prime, x^{\prime 0}} f)(x) \ \text{ с } \ \delta^\prime = \delta^{-1},
x^{\prime 0} =-\delta^{-1} x^0.
\end{equation*}

Отметим, что при $ 1 \le p \le \infty $ для $ f \in L_p(x^0 +\delta D), $ где
$ D $ -- область в $ \R^d, \delta \in \R_+^d, x^0 \in \R^d, $ имеет место
равенство
\begin{equation*} \tag{1.4.1}
\| h_{\delta,x^0} f\|_{L_p(D)} = \delta^{-p^{-1} \e} \|f\|_{L_p(x^0 +\delta D)},
\end{equation*}
а, следовательно, для $ f \in L_p(D) $ выполняется равенство
\begin{equation*} \tag{1.4.2}
\| h_{\delta,x^0}^{-1} f\|_{L_p(x^0 +\delta D)} = \delta^{p^{-1} \e} \|f\|_{L_p(D)}.
\end{equation*}

Лемма 1.4.1

Пусть $ d \in \N, D $ -- область в $ \R^d, \alpha \in \R_+^d, 1 \le p < \infty,
1 \le \theta \le \infty, \delta \in \R_+^d, x^0 \in \R^d. $ Тогда существуют
константы $ c_1(d,\alpha,p,\delta) > 0, c_2(d,\alpha,p,\delta) > 0 $
такие, что для любой функции $ f \in (B_{p,\theta}^\alpha)^\prime(x^0 +\delta D) $
соблюдается неравенство
\begin{equation*} \tag{1.4.3}
\| h_{\delta, x^0} f \|_{(B_{p,\theta}^\alpha)^\prime(D)} \le
c_1 \| f \|_{(B_{p,\theta}^\alpha)^\prime(x^0 +\delta D)},
\end{equation*}
а для $ f \in (B_{p,\theta}^\alpha)^\prime(D) $ выполняется неравенство
\begin{equation*} \tag{1.4.4}
\| h_{\delta, x^0}^{-1} f \|_{(B_{p,\theta}^\alpha)^\prime(x^0 +\delta D)} \le
c_2 \| f \|_{(B_{p,\theta}^\alpha)^\prime(D)}.
\end{equation*}

Определение области $\alpha$-типа см. в [14], в [6] (см. также [14]) установлена
следующая теорема.

Теорема 1.4.2

Пусть $ d \in \N, \alpha \in \R_+^d, 1 \le p < \infty, 1 \le \theta \le \infty $
и $ D \subset \R^d $ -- ограниченная область $ \alpha $-типа.
Тогда существует непрерывное линейное отображение
$ E^{d, \alpha, p, \theta,D}: (B_{p, \theta}^\alpha)^\prime(D) \mapsto
L_p(\R^d) $ такое, что для любой функции $ f \in (B_{p, \theta}^\alpha)^\prime(D) $
соблюдается равенство
\begin{equation*} \tag{1.4.5}
(E^{d, \alpha, p, \theta,D} f) \mid_{D} = f,
\end{equation*}
а при $ \l \in \Z_+^d: \l < \alpha, $ существует константа
$ c_3(d,\alpha,p,\theta,D,\l) > 0, $ для которой выполняется неравенство
\begin{equation*} \tag{1.4.6}
\| E^{d, \alpha, p, \theta,D} f\|_{(B_{p, \theta}^\alpha)^{\l}(\R^d)} \le
c_3 \| f\|_{(B_{p, \theta}^\alpha)^\prime(D)}.
\end{equation*}

Следствие

В условиях теоремы 1.4.2 для любой функции $ f \in (B_{p, \theta}^\alpha)^\prime(D) $
выполняется неравенство
\begin{equation*} \tag{1.4.7}
\| E^{d, \alpha, p, \theta,D} f\|_{(B_{p, \theta}^\alpha)^\prime(\R^d)} \le
c_3 \| f\|_{(B_{p, \theta}^\alpha)^\prime(D)}.
\end{equation*}

Для получения (1.4.7) достаточно с учётом (1.1.3) соединить (1.1.4) и
(1.4.6).
\bigskip

\centerline{\S 2. Оценка бернштейновского $n$-поперечgика шара}
\centerline{$ B((B_{p,\theta}^\alpha)^\prime(D)) $ в пространстве $L_q(D) $}
\centerline{в ограниченной области $ D \ \alpha$-типа}
\bigskip

2.1. В этом пункте получим оценку сверху интересующей нас величины.

Теорема 2.1.1

Пусть $ d \in \N, \alpha \in \R_+^d, 1 \le p < \infty, 1 \le \theta \le \infty,
1 \le q \le \infty $ удовлетворяют условию (1.2.10) и $ D \subset \R^d $ --
ограниченная область $\alpha$-типа. Тогда существуют константы
$ c_1(d,\alpha,p,\theta,q,D) > 0 $ и $ n_0(d,\alpha) > 0 $ такие, что при
$ n \ge n_0 $ для $ C = B((B_{p,\theta}^\alpha)^\prime(D)), X = L_q(D) $
соблюдаются неравенства
\begin{multline*} \tag{2.1.1}
b_n(C,X) \le c_1 \cdot
\begin{cases}
n^{-1 /(\alpha^{-1}, \e)}, \parbox[t]{7cm}{ при $ 1 \le q \le p \le 2$ или $ 1 \le q = p <
\infty$ или ( $ 1 \le p < q \le \infty$ и $ 1 -(\alpha^{-1}, \e) /p +
(\alpha^{-1}, \e) /q >0 $ );} \\
n^{-\!(1 /(\alpha^{-1}, \e) -1 /p +1/2)},  \parbox[t]{5.5cm}{ при $ 1 \le q \le 2 < p < \infty $
и $ 1 -(\alpha^{-1}, \e) /p > 0$;} \\
n^{-\!(1 /(\alpha^{-1}, \e) -1 /p +1 /q)},  \parbox[t]{5.5cm}{ при $ 2 \le q < p < \infty $
и $ 1 -(\alpha^{-1}, \e) /p +(\alpha^{-1}, \e) /q -(\alpha^{-1}, \e) /2 >0.$}
\end{cases}
\end{multline*}

Доказательство.

Сначала проведём оценку сверху $n$-поперечника по Бернштейну для $ C =
B((H_p^\alpha)^\prime(I^d)) $ в пространстве $ X = L_q(I^d)$. Эта оценка
осуществляется с помощью лемм 2.1.2--2.1.4, теоремы 1.3.3 и других средств.

Лемма 2.1.2

Пусть выполнены условия леммы 1.2.6, а также пусть
$ C = B((H_p^\alpha)^\prime(I^d)),\ X = L_q(I^d). $ Тогда для любого $ n \in \N$
существует $ m_n \in \N$ такое, что при любом $ k \ge m_n$ справедлива оценка
\begin{equation*} \tag{2.1.2}
b_n(C,X) \le 2 b_n(E_k^{d,l -\e,\alpha}(C), \mathcal P_k^{d,l -\e,\alpha} \bigcap X),
\end{equation*}
где $ \mathcal P_k^{d,l -\e,\alpha} \bigcap X$, как обычно, означает, что в
пространстве $ \mathcal P_k^{d,l -\e,\alpha} $ рассматривается норма,
индуцированная из пространства $X$.

Доказательство.

Прежде всего, заметим, что в условиях леммы ввиду (1.2.10) последовательность
$$
c_0 2^{-k(1 -(\alpha^{-1},\e)(p^{-1} -q^{-1})_+)} \to 0 \text{ при $ k \to +\infty$,}
$$
где $ c_0 $ --- это константа $ c_8 $ из (1.2.11).
Для $ n \in \N$ выберем число $ m_n \in \N$ так, чтобы при $ k \ge m_n$
соблюдалось неравенство
\begin{equation*} \tag{2.1.3}
c_0 2^{-k(1 -(\alpha^{-1},\e)(p^{-1} -q^{-1})_+ )} \le \frac{1}{2} b_n(C,X).
\end{equation*}

Фиксировав $ n \in \N$, возьмём для произвольного $ \delta:
0 < \delta < \frac{1}{2} b_n(C,X),$ подпространство $ L \in \mathcal M_n(X)$
и число $ \epsilon > b_n(C,X) -\delta > \frac{1}{2} b_n(C,X)$ такие, что
$ (\epsilon \cdot B(X)) \bigcap L \subset C.$

Тогда при $ k \ge m_n$ для любого элемента $ f \in L$ в силу сказанного с
учётом (1.2.11) и (2.1.3) имеем
\begin{multline*}
\left\| E_k^{d,l -\e,\alpha} \biggl(\frac{\epsilon} {\|f\|_X} f\biggr)\right\|_X
 \ge \left\| \frac{\epsilon} {\|f\|_X} f\right\|_X -\left\| \frac{\epsilon} {\|f\|_X} f -
E_k^{d,l -\e,\alpha} \biggl(\frac{\epsilon} {\|f\|_X} f\biggr)\right\|_X \\
 \ge \epsilon -c_0 2^{-k(1 -(\alpha^{-1},\e)(p^{-1} -q^{-1})_+)} \ge \epsilon
-\frac{1}{2} b_n(C,X)
\end{multline*}
или
\begin{equation*} \tag{2.1.4}
\bigl\| E_k^{d,l -\e,\alpha} f\bigr\|_X \ge \biggl(1 -\frac{b_n(C,X)}
{2 \epsilon}\biggr) \|f\|_X, \ f \in L.
\end{equation*}

Принимая во внимание, что $ 1 -\frac{b_n(C,X)} {2 \epsilon} >0$, из
(2.1.4) заключаем, что для $ f \in L \setminus \{0\} $ его образ
$ E_k^{d,l -\e,\alpha} f \ne 0$. Поэтому
$ E_k^{d,l -\e,\alpha}(L) \in \mathcal M_n(\mathcal P_k^{d,l -\e,\alpha})$.

Кроме того, если $ g = E_k^{d,l -\e,\alpha} f$, где $ f \in L$ и
$ \|g\|_X \le \epsilon -\frac{1}{2} b_n(C,X)$, то в силу (2.1.4) норма
$ \|f\|_X \le \epsilon $ и, следовательно, $ f \in C$, т.е.
\begin{equation*}
(\epsilon -\frac{1}{2} b_n(C,X)) \cdot B(X) \bigcap E_k^{d,l -\e,\alpha}(L)
\subset E_k^{d,l -\e,\alpha}(C).
\end{equation*}

Поэтому
\begin{multline*}
b_n( E_k^{d,l -\e,\alpha}(C), \mathcal P_k^{d,l -\e,\alpha} \cap X)
\ge \epsilon -\frac{1}{2} b_n(C,X) > \\
b_n(C,X) -\delta 
-\frac{1}{2} b_n(C,X) = \frac{1}{2} b_n(C,X) -\delta.
\end{multline*}
Откуда, ввиду произвольности $ \delta >0,$ вытекает (2.1.2). $\square$

Лемма 2.1.3

В условиях леммы 2.1.2 для любых $ k \in\Z_+,\ j \in \N$ при $ R_k^{d,l -\e,\alpha}
\le n \le \frac {1}{2} R_{k +j}^{d,l -\e,\alpha}$ имеет место оценка
\begin{equation*} \tag{2.1.5}
b_{2n}(E_{k +j}^{d,l -\e,\alpha}(C),\mathcal P_{k +j}^{d,l -\e,\alpha}
\bigcap X) \le b_n(\mathfrak E_{k,k +j}^{d,l -\e,\alpha}(C),
\mathfrak P_{k,k +j}^{\prime \,d,l -\e,\alpha} \bigcap X).
\end{equation*}

Доказательство.

Пусть $ L \in \mathcal M_{2n}(\mathcal P_{k +j}^{d,l -\e,\alpha})$ и $ \epsilon
\ge 0$ таковы, что $ (\epsilon \cdot B(X)) \bigcap L \subset
E_{k +j}^{d,l -\e,\alpha}(C)$. Поскольку
$$
\dim \mathcal P_{k +j}^{d,l -\e,\alpha} \ge
\dim \mathfrak P_{k,k +j}^{\prime \,d,l -\e,\alpha} +\dim L -
\dim \left(L \bigcap \mathfrak P_{k,k +j}^{\prime \,d,l -\e,\alpha} \right),
$$
то (1.2.15) даёт
\begin{equation*}
\dim \mathcal P_{k +j}^{d,l -\e,\alpha} \ge \dim \mathcal P_{k +j}^{d,l -\e,\alpha}
-\dim \mathcal P_k^{d,l -\e,\alpha} +\dim L -\dim \left(L \bigcap
\mathfrak P_{k,k +j}^{\prime \,d,l -\e,\alpha} \right),
\end{equation*}
т.е.
$$
\dim \left(L \bigcap \mathfrak P_{k,k +j}^{\prime \,d,l -\e,\alpha} \right) \ge
\dim L -\dim \mathcal P_k^{d,l -\e,\alpha} = 2n -R_k^{d,l -\e,\alpha} \ge n.
$$

Таким образом, существует подпространство $ M \in
\mathcal M_n(\mathfrak P_{k,k +j}^{\prime \,d,l -\e,\alpha} )$ такое, что
$ M \subset L$. Но тогда
$ (\epsilon \cdot B(X)) \bigcap M \subset (\epsilon \cdot B(X)) \bigcap L
\subset E_{k +j}^{d,l -\e,\alpha}(C)$.
Поэтому для $ f \in (\epsilon \cdot B(X)) \bigcap M \subset
\mathfrak P_{k,k +j}^{\prime \,d,l -\e,\alpha} )$ имеем
$ f = E_{k +j}^{d,l -\e,\alpha} g,\ g \in C$, и с учётом (1.2.14) $ f =
\mathfrak E_{k,k +j}^{d,l -\e,\alpha} f $.
Следовательно, пользуясь (1.2.12), получаем, что
\begin{multline*}
f  = \mathfrak E_{k,k +j}^{d,l -\e,\alpha} E_{k +j}^{d,l -\e,\alpha} g
= ((E_{k +j}^{d,l -\e,\alpha})^2 -E_k^{d,l -\e,\alpha} E_{k +j}^{d,l -\e,\alpha}) g = \\
(E_{k +j}^{d,l -\e,\alpha} -E_k^{d,l -\e,\alpha}) g 
 = \mathfrak E_{k,k +j}^{d,l -\e,\alpha} g,\ \text{ т.е. } 
(\epsilon \cdot B(X)) \bigcap M \subset
\mathfrak E_{k,k +j}^{d,l -\e,\alpha}(C).
\end{multline*}

А это значит, что $ \epsilon \le b_n\bigl(\mathfrak E_{k,k +j}^{d,l -\e,\alpha}(C),
\mathfrak P_{k,k +j}^{\prime \,d,l -\e,\alpha} \bigcap X\bigr)$.

Отсюда в силу произвольности $ L $ и $ \epsilon $ вытекает (2.1.5). $\square$

Лемма 2.1.4

Пусть $ d \in \N, \alpha \in \R_+^d, l = l(\alpha) \in \N^d, 1 \le p < \infty,
C = B((H_p^\alpha)^\prime(I^d)). $ Тогда существует константа $ c_2(d,\alpha,p) >0$
такая, что для любых $ k \in \Z_+,\ j_0 \in \ N$ имеет место включение
\begin{equation*} \tag{2.1.6}
\mathfrak E_{k,k +j_0}^{d,l -\e,\alpha}(C) \subset \bigcap_{j =1}^{j_0}
\left\{\,f \in \mathfrak P_{k,k +j_0}^{\prime \,d,l -\e,\alpha}:
\bigl\| \mathcal E_{k +j}^{d,l -\e,\alpha} f\bigr\|_{L_p(I^d)} \le c_2
2^{-(k +j)} \,\right\}.
\end{equation*}

Доказательство.

В условиях леммы при $ k \in \Z_+,\ j_0 \in \ N $ пусть
$ f = \mathfrak E_{k,k +j_0}^{d,l -\e,\alpha} g$, где $ g \in C$. Тогда
$ f \in \mathfrak P_{k,k +j_0}^{\prime \,d,l -\e,\alpha} $ и
$ f = \sum\limits_{j =1}^{j_0} \mathcal E_{k +j}^{d,l -\e,\alpha} g, $ а в силу
(1.2.13) для $ i =1,\ldots,j_0$ имеем
$$
\mathcal E_{k +i}^{d,l -\e,\alpha} f = \sum_{j =1}^{j_0} \mathcal E_{k +i}^{d,l -\e,\alpha}
\mathcal E_{k +j}^{d,l -\e,\alpha} g = \mathcal E_{k +i}^{d,l -\e,\alpha} g.
$$
Поэтому, пользуясь (1.2.11), получаем
\begin{multline*}
\bigl\|\mathcal E_{k +j}^{d,l -\e,\alpha} f\bigr\|_{L_p(I^d)}
= \bigl\| \mathcal E_{k +j}^{d,l -\e,\alpha} g\bigr\|_{L_p(I^d)} \le \bigl\| g
-E_{k +j}^{d,l -\e,\alpha} g\bigr\|_{L_p(I^d)} +\bigl\| g -
E_{k +j -1}^{d,l -\e,\alpha} g\bigr\|_{L_p(I^d)} \\
\le c_3 2^{-(k +j)} +c_3 2^{-(k +j -1)}
 \le c_2 2^{-(k +j)},\ j =1,\ldots,j_0. \square
\end{multline*}

Продолжим доказательство теоремы 2.1.1.
В условиях леммы 2.1.2 пусть $ n \in \N: n \ge n_0 = 2 R_0^{d,l -\e,\alpha}$.
Фиксировав $ k \in \Z_+$, для которого выполняется неравенство
\begin{equation*} \tag{2.1.7}
2 R_k^{d,l -\e,\alpha} = 2 R_k \le n < 2 R_{k +1},
\end{equation*} (здесь и ниже для упрощения записи при проведении выкладок
верхние индексы $ d, l -\e, \alpha$ опускаем),
принимая во внимание лемму 2.1.2 и (1.2.5), выберем $ j_n \in \N $ так,
чтобы было $ j_n \ge m_{2 R_k}$ и $ R_k \le \frac {1}{2} R_{k +j_n}$.

Тогда, благодаря (2.1.7), (1.3.4), (2.1.2), (2.1.5), получаем
\begin{multline*} \tag{2.1.8}
b_n(C,X) \le b_{2 R_k}(C,X) \le 2 b_{2 R_k}(E_{k +j_n}(C),
\mathcal P_{k +j_n} \bigcap X) \\
\le 2 b_{R_k}(\mathfrak E_{k,k +j_n}(C), \mathfrak P_{k, k +j_n}^{\prime}
\bigcap X).
\end{multline*}

Далее, обозначим через $ C^\prime,\ C_j^\prime,\ j =1,\ldots,j_n$, множества
\begin{gather*}
C^\prime = \bigcap_{j =1}^{j_n} C_j^\prime, \\
C_j^\prime = \left\{\,f \in \mathfrak P_{k,k +j_n}^{\prime}: \|\mathcal E_{k +j} f\|_{L_p(I^d)} \le
c_2 2^{-(k +j)} \,\right\} \text{ (см. (2.1.6))}.
\end{gather*}
Тогда для $ f \in \mathfrak P_{k,k +j_n}^{\prime}$ имеем
\begin{multline*}
\mu_{C_j^\prime}(f)
 = \inf\left\{\,t >0: \bigl\|\mathcal E_{k +j}(t^{-1} f)\bigr\|_{L_p(I^d)}
\le c_2 2^{-(k +j)}\,\right\} \\
= \inf\biggl\{\,t >0: \frac{\bigl\|\mathcal E_{k +j} f\bigr\|_{L_p(I^d)}}
{c_2 2^{-(k +j)}} \le t\,\biggr\}
= \frac{\bigl\|\mathcal E_{k +j} f\bigr\|_{L_p(I^d)}}
{c_2 2^{-(k +j)}}.
\end{multline*}

Из (2.1.8), принимая во внимание сказанное и учитывая (2.1.6), (1.3.2),
(1.3.7), (1.3.6), выводим
\begin{multline*} \tag{2.1.9}
b_n(C,X)
 \le 2 b_{R_k}(C^\prime, \mathfrak P_{k,k +j_n}^{\prime} \bigcap X) =
2 \sup_{L \in \mathcal M_{R_k}\left(\mathfrak P_{k,k +j_n}^{\prime}\right)}
\inf_{f \in L \setminus \{0\}} \frac { \|f\|_X } { \mu_{C^\prime}(f)} \\
 = 2 \sup_{L \in \mathcal M_{R_k}\left(\mathfrak P_{k,k +j_n}^{\prime}\right)}
\inf_{f \in L \setminus \{0\}} \|f\|_X \cdot \left(\max_{j =1,\ldots,j_n}
\mu_{C_j^\prime}(f)\right)^{\!-1} \\
 = 2 \sup_{L \in \mathcal M_{R_k}\left(\mathfrak P_{k,k +j_n}^{\prime}\right)}
\inf_{f \in L \setminus \{0\}} \|f\|_{L_q(I^d)}  \left(\max_{j =1,\ldots,j_n}
\frac {\bigl\|\mathcal E_{k +j} f\bigr\|_{L_p(I^d)}}
{c_2 2^{-(k +j)}} \right)^{\!-1}.
\end{multline*}

При проведении дальнейших рассуждений рассмотрим несколько случаев
соотношений между $p$ и $q$.

Случай 1. $ q \le p \le 2$.

В этом случае, используя (1.2.11), имеем
\begin{multline*}
b_{R_k}(\mathfrak E_{k,k +j_n}(C), \mathfrak P_{k,k +j_n}^{\prime} \bigcap X)
\le \sup_{f \in C} \| \mathfrak E_{k,k +j_n} f\|_X \\
 \le \sup_{f \in C} \biggl(\bigl\| f -E_{k +j_n} f\bigr\|_{L_q(I^d)} +
\bigl\| f -E_k f\bigr\|_{L_q(I^d)} \biggr) \\
 \le c_4 2^{-(k +j_n)} +c_4 2^{-k} \le 2 c_4 2^{-k}.
\end{multline*}
Соединяя эту оценку с (2.1.8) и учитывая (1.2.5), (2.1.7), получаем
(2.1.1) при $ q \le p \le 2$.

Случай $1^\prime$. $ p = q$.

В этом случае доказательство дословно повторяет доказательство в случае 1.

Случай 2. $ p > \max(2,q)$.

В этом случае, учитывая (1.2.14), (1.2.4) и (1.2.6), для $ f \in
\mathfrak P_{k,k +j_n}^{\prime} $ имеем
\begin{multline*} \tag{2.1.10}
\|f\|_{L_q(I^d)}
 = \| \mathfrak E_{k,k +j_n} f\|_{L_q(I^d)} = \biggl\| \sum_{j =1}^{j_n}
\mathcal E_{k +j} f\biggr\|_{L_q(I^d)} \\
 \le \sum_{j =1}^{j_n} \| \mathcal E_{k +j} f\|_{L_q(I^d)} \le \sum_{j =1}^{j_n}
c_5 2^{-(k +j)(\alpha^{-1}, \e) /q} \cdot \bigl\| \I_{k +j}
\mathcal E_{k +j} f\bigr\|_{l_q^{R_{k +j}}}.
\end{multline*}

Из (2.1.10), применяя неравенство (1.1.1), затем (1.2.5) и неравенство
Гёльдера, получаем
\begin{multline*} \tag{2.1.11}
\|f\|_{L_q(I^d)}
 \le c_5 \sum_{j =1}^{j_n} 2^{-(k +j)(\alpha^{-1}, \e) /q}
\bigl(R_{k +j}\bigr)^{(1/q -1/2)_+} \cdot \| \I_{k +j} \mathcal E_{k +j} f\|_{l_2^{R_{k +j}}} \\
 \le c_6 \sum_{j =1}^{j_n}
2^{-(k +j)(\alpha^{-1}, \e) /q +(k +j)(\alpha^{-1}, \e)(1/q -1/2)_+ +\epsilon j} \cdot
\| \I_{k +j} \mathcal E_{k +j} f\|_{l_2^{R_{k +j}}} \cdot 2^{-\epsilon j}\\
 \le c_6 \left(\sum_{j =1}^{j_n} \biggl(2^{-(k +j)(\alpha^{-1}, \e) /q +
(k +j)(\alpha^{-1}, \e)(1/q -1/2)_+ +\epsilon j}
\cdot \| \I_{k+j} \mathcal E_{k +j} f\|_{l_2^{R_{k +j}}}\biggr)^2\right)^{1/2} \\
 \times \biggl(\sum_{j =1}^{j_n} 2^{-2 \epsilon j} \biggr)^{1/2}
\le c_7 \left(\sum_{j =1}^{j_n} \biggl(2^{-(k +j)(\alpha^{-1}, \e) /q +
(k +j)(\alpha^{-1}, \e)(1/q -1/2)_+ +\epsilon j}
\right.\\
 \times \left.\| \I_{k +j} \mathcal E_{k +j} f\|_{l_2^{R_{k +j}}}\biggr)^2\right)^{1/2},
\end{multline*}
где $ \epsilon >0 $ будет выбрано ниже.

Пользуясь (1.2.6), (1.1.1), находим для $ f \in \mathfrak P_{k,k +j_n}^{\prime} $
\begin{multline*} \tag{2.1.12}
\max_{j =1,\ldots,j_n}
 \frac{\bigl\| \mathcal E_{k +j} f\bigr\|_{L_p(I^d)}}
{c_2 2^{-(k +j)}} \ge \max_{j =1,\ldots,j_n}
\frac{c_8 2^{-(k +j)(\alpha^{-1}, \e) /p} \bigl\| \I_{k +j}
\mathcal E_{k +j} f\bigr\|_{l_p^{R_{k +j}}}} {c_2 2^{-(k +j)}} \\
 \ge c_9 \max_{j =1,\ldots,j_n}
\frac{ \bigl\| \I_{k +j} \mathcal E_{k +j} f\bigr\|_{l_\infty^{R_{k +j}}}}
{2^{-(k +j)(1 -(\alpha^{-1}, \e) /p)}}.
\end{multline*}

Определим линейный оператор $ \mathcal A^n: \mathfrak P_{k,k +j_n}^{\prime}
\mapsto \R^{\mathcal R_n} $, полагая для $ f \in \mathfrak P_{k,k +j_n}^{\prime}$
его образ
\begin{gather*}
\mathcal A^n f = y = \{\,y_j = \{\,y_j^{\nu,\lambda} \in \R: \nu \in
\Nu_{0,2^{\kappa(k +j,\alpha)} -\e}^d,\ \lambda \in \Z_+^d(l -\e)\,\}, \
j =1,\ldots,j_n\,\}, \\
\text{ где } y_j = 2^{-(k +j)(\alpha^{-1}, \e) /q +(k +j)(\alpha^{-1}, \e)(1/q -1/2)_+ +\epsilon j} \cdot \I_{k +j}
\mathcal E_{k +j} f \in \R^{R_{k +j}}, \\
\mathcal R_n = \sum_{j =1}^{j_n} R_{k +j}.
\end{gather*}

Заметим, что ядро $ \Ker \mathcal A^n = \{0\}$, ибо, если для $ f \in
\mathfrak P_{k,k +j_n}^{\prime} $ его образ $ \mathcal A^n f =0$, т.е.
$ \I_{k +j} \mathcal E_{k +j} f =0$ для $ j =1,\ldots,j_n$ , то
$ \mathcal E_{k +j} f =0$ для $ j =1,\ldots,j_n$ (см. лемму 1.2.2) и, значит,
(см. (1.2.14)) $ f = \sum\limits_{j =1}^{j_n}
\mathcal E_{k +j} f =0$.

Из (2.1.11) и (2.1.12) с учётом введённых обозначений для $ f \in
\mathfrak P_{k,k +j_n}^{\prime}$ соответственно получаем
\begin{equation*} \tag{2.1.13}
\|f\|_{L_q(I^d) }\le c_7 \biggl(\sum_{j =1}^{j_n}
\biggl(\|y_j\|_{l_2^{R_{k +j}}}\biggr)^{\!2}\biggr)^{\!1/2} = c_7
\|y\|_{l_2^{\mathcal R_n}} = c_7 \| \mathcal A^n f\|_{l_2^{\mathcal R_n}},
\end{equation*}

\begin{multline*} \tag{2.1.14}
\max_{j =1,\ldots,j_n}
\frac{ \bigl\| \mathcal E_{k +j} f\bigr\|_{L_p(I^d)} }
{ c_2 2^{-(k +j)} } \\
 \ge c_9 \max_{j =1,\ldots,j_n}
\frac{ \|y_j\|_{l_\infty^{ R_{k +j} } } }
{ 2^{-(k+j)(1 -(\alpha^{-1}, \e) /p +(\alpha^{-1}, \e) /q -(\alpha^{-1}, \e)
(1 /q -1/2)_+) +\epsilon j} } \\
 = c_9 \max_{ \{\,(j,\nu,\lambda): j =1,\ldots,j_n,\ \nu \in
\Nu_{0, 2^{\kappa(k +j, \alpha)} -\e}^d,\ \lambda \in \Z_+^d(l -\e)\,\} }
\frac{ | y_j^{\nu,\lambda}| } {\rho_j^{\nu,\lambda}} \\
 = c_9 \mu_{\mathcal B}(y) = c_9 \mu_{\mathcal B}(\mathcal A^n f),
\end{multline*}
где
\begin{multline*}
\mathcal B = \{\,y \in \R^{\mathcal R_n}: | y_j^{\nu,\lambda}| \le \rho_j^{\nu,\lambda} =
\rho_j = 2^{-(k +j)(1 -(\alpha^{-1}, \e) /p +(\alpha^{-1}, \e) /q -
(\alpha^{-1}, \e)(1 /q -1/2)_+) +\epsilon j}, \\
j =1,\ldots,j_n,\ \nu \in \Nu_{0,2^{\kappa(k +j, \alpha)} -\e}^d,\ \lambda
\in \Z_+^d(l -\e)\,\}.
\end{multline*}

Соединяя (2.1.9), (2.1.13) и (2.1.14) и учитывая, что для $ L \in
\mathcal M_{R_k} \bigl(\mathfrak P_{k,k +j_n}^{\prime}\bigr) $ его образ
$ \mathcal A^n(L) \in \mathcal M_{R_k}(\R^{\mathcal R_n})$, а также, принимая во
внимание (1.3.7), находим, что
\begin{multline*} \tag{2.1.15}
b_n(C,X)
 \le 2 \sup_{L \in \mathcal M_{R_k}\bigl(\mathfrak P_{k,k +j_n}^{\prime}\bigr)}
\inf_{f \in L \setminus \{0\}} \frac{ c_7 \| \mathcal A^n f \|_{l_2^{\mathcal R_n}} }
{c_9 \mu_{\mathcal B}(\mathcal A^n f)} \\
 = c_{10} \sup_{L \in \mathcal M_{R_k}\bigl(\mathfrak P_{k,k +j_n}^{\prime}\bigr)}
\inf_{y \in \mathcal A^n(L) \setminus \{0\}} \frac{ \|y\|_{l_2^{\mathcal R_n}} }
{ \mu_{\mathcal B}(y)} \\
 = c_{10} \sup_{ \left\{ M = \mathcal A^n(L): L \in
\mathcal M_{R_k}\bigl(\mathfrak P_{k,k +j_n}^{\prime}\bigr)\right\}}
\inf_{ y \in M \setminus \{0\}} \frac{ \|y\|_{l_2^{\mathcal R_n}} }
{ \mu_{\mathcal B}(y)} \\
 \le c_{10} \sup_{M \in \mathcal M_{R_k}(\R^{\mathcal R_n})}
\inf_{y \in M \setminus \{0\}} \frac{ \|y\|_{l_2^{\mathcal R_n}} }
{ \mu_{\mathcal B}(y)} 
 = c_{10} b_{R_k}(\mathcal B, l_2^{\mathcal R_n}).
\end{multline*}

Применяя (1.3.8), (1.2.5), выводим
\begin{multline*} \tag{2.1.16}
b_{R_k}(\mathcal B, l_2^{\mathcal R_n})
 \le \biggl(\frac{2}{R_k} \sum_{\{\,j = 1,\ldots,j_n,\ \nu \in
\Nu_{0, 2^{\kappa(k +j,\alpha)} -\e}^d,\ \lambda \in \Z_+^d(l -\e)\,\}}
(\rho_j^{\nu,\lambda})^2 \biggr)^{1/2} \\
 = \biggl( \frac{2}{R_k} \sum_{j =1}^{j_n} \sum_{\{\,\nu \in
\Nu_{0,2^{\kappa(k +j,\alpha)} -\e}^d,\ \lambda \in \Z_+^d(l -\e)\,\}}
(\rho_j^{\nu,\lambda})^2 \biggr)^{1/2} = \\
\biggl( \frac{2}{R_k} \sum_{j =1}^{j_n} R_{k +j} (\rho_j)^2 \biggr)^{1/2} 
 \le c_{11} \biggl(\sum_{j =1}^{j_n} 2^{j(\alpha^{-1}, \e)} (\rho_j)^2 \biggr)^{1/2} = \\
c_{11} \biggl(\sum_{j =1}^{j_n}
2^{-2(k +j)(1 -(\alpha^{-1}, \e) /p +(\alpha^{-1}, \e) /q -(\alpha^{-1}, \e)
(1 /q -1/2)_+) +2 \epsilon j +(\alpha^{-1}, \e) j} \biggr)^{1/2} \le \\
c_{11} 2^{-k(1 -(\alpha^{-1}, \e) /p +(\alpha^{-1}, \e) /q -(\alpha^{-1}, \e)
(1 /q -1/2)_+)} \times \\
\biggl(\sum_{j =1}^{\infty}
2^{-2j(1 -(\alpha^{-1}, \e) /p +(\alpha^{-1}, \e) /q -(\alpha^{-1}, \e)
(1 /q -1/2)_+ -(\alpha^{-1}, \e) /2 -\epsilon)} \biggr)^{1/2} = \\
c_{11} 2^{-k(1 -(\alpha^{-1}, \e) /p +(\alpha^{-1}, \e) /q -(\alpha^{-1}, \e)
(1 /q -1/2)_+)} \times\\
\biggl(\sum_{j =1}^{\infty}
2^{-2j(1 -(\alpha^{-1}, \e) /p -(\alpha^{-1}, \e)(1 /2 -1 /q)_+ -\epsilon)} \biggr)^{1/2} = \\
= c_{12} 2^{-k(1 -(\alpha^{-1}, \e) /p +(\alpha^{-1}, \e) /q -(\alpha^{-1}, \e)
(1 /q -1/2)_+)},
\end{multline*}
где $ 0 < c_{12} < \infty$, а $ \epsilon > 0$ взято таким, чтобы было
$ 1 -(\alpha^{-1}, \e) /p -(\alpha^{-1}, \e)(1 /2 -1 /q)_+ -\epsilon >0$, что
возможно при условии $ 1 -(\alpha^{-1}, \e) /p -(\alpha^{-1}, \e)(1 /2 -1 /q)_+
> 0.$

Подставляя оценку (2.1.16) в (2.1.15) и учитывая (1.2.5), (2.1.7), получаем
(2.1.1) при $ p > \max(2,q) $ и условии $ 1 -(\alpha^{-1}, \e) /p -
(\alpha^{-1}, \e)(1 /2 -1 /q)_+ > 0. $

Случай 3. $ p < q. $

Из (2.1.10)  с помощью неравенства Гёльдера для $ f \in
\mathfrak P_{k,k +j_n}^{\prime} $ получаем
\begin{multline*} \tag{2.1.17}
\|f\|_{L_q(I^d)}  \le c_5 \sum_{j =1}^{j_n} 2^{-(k +j)(\alpha^{-1}, \e) /q
+\epsilon j} \| \I_{k +j} \mathcal E_{k +j} f\|_{l_q^{R_{k +j}}} 2^{-\epsilon j} \\
 \le c_5 \left(\sum_{j =1}^{j_n} \biggl(2^{-(k +j)(\alpha^{-1}, \e) /q +
\epsilon j} \| \I_{k +j} \mathcal E_{k +j} f\|_{l_q^{R_{k +j}}}\biggr)^q \right)^{1/q}
\biggl(\sum_{j =1}^{j_n} 2^{-\epsilon j q^\prime} \biggr)^{1 /q^\prime} \\
 \le c_{13} \left(\sum_{j =1}^{j_n} \biggl(2^{-(k +j)(\alpha^{-1}, \e) /q +
\epsilon j} \| \I_{k +j} \mathcal E_{k +j} f\|_{l_q^{R_{k +j}}}\biggr)^q \right)^{1/q},
\end{multline*}
где $ \epsilon > 0$ будет выбрано ниже.

Далее, ясно, что для $ f \in \mathfrak P_{k,k +j_n}^{\prime} $ при $ \delta >0$,
которое будет указано ниже, соблюдается неравенство
\begin{multline*} \tag{2.1.18}
 \left(\sum_{j =1}^{j_n} \biggl(\frac{ \| \mathcal E_{k +j} f\|_{L_p(I^d)} }
{ c_2 2^{-(k +j) +\delta j}} \biggr)^p \right)^{1 /p} 
 \le \biggl(\sum_{j =1}^{j_n} 2^{-\delta j p} \biggr)^{1 /p}
\cdot \max_{j =1,\ldots,j_n} \frac{ \| \mathcal E_{k +j} f\|_{L_p(I^d)} }
{ c_2 2^{-(k +j)} } \\
 \le c_{14} \max_{j =1,\ldots,j_n} \frac{ \| \mathcal E_{k +j} f\|_{L_p(I^d)} }
{ c_2 2^{-(k +j)} }.
\end{multline*}

Из (2.1.18) с использованием (1.2.6) выводим
\begin{equation*} \tag{2.1.19}
\max_{j =1,\ldots,j_n} \frac{ \| \mathcal E_{k +j} f\|_{L_p(I^d)} }
{ c_2 2^{-(k +j)} } 
\ge \frac{1}{c_{14}} \left(\sum_{j =1}^{j_n} \left( \frac{c_{15}
2^{-(k +j)(\alpha^{-1}, \e) /p} \| \I_{k +j} \mathcal E_{k +j} f\|_{l_p^{R_{k +j}}} }
{ c_2 2^{-(k +j) +\delta j}} \right)^p \right)^{1 /p}.
\end{equation*}

Определим линейный оператор $ \mathcal A^n: \mathfrak P_{k,k +j_n}^{\prime}
\mapsto \R^{\mathcal R_n},\ \mathcal R_n = \sum\limits_{j =1}^{j_n} R_{k +j}$,
полагая для $ f \in \mathfrak P_{k,k +j_n}^{\prime}$ его образ равным
\begin{multline*}
\mathcal A^n f = y = \{\,y_j = \{\,y_j^{\nu,\lambda} \in \R,\ \nu \in
\Nu_{0, 2^{\kappa(k +j, \alpha)} -\e}^d,\
\lambda \in \Z_+^d(l -\e)\,\},\ j =1,\ldots,j_n\,\}, \\
\text{ где } y_j = 2^{-(k +j)(\alpha^{-1}, \e) /q +\epsilon j} \cdot \I_{k +j}
\mathcal E_{k +j} f,\ j =1,\ldots,j_n.
\end{multline*}

Как и в случае 2, если для $ f \in \mathfrak P_{k,k +j_n}^{\prime} $ его
образ $ \mathcal A^n f = 0, $ то $ f =0$.

Пользуясь этими обозначениями, из (2.1.17) и (2.1.19) для $ f \in
\mathfrak P_{k,k +j_n}^{\prime} $ соответственно имеем
\begin{equation*} \tag{2.1.20}
\|f\|_{L_q(I^d) } \le c_{13} \biggl(\sum_{j =1}^{j_n}
\biggl( \|y_j\|_{l_q^{R_{k +j}}} \biggr)^q \biggr)^{1 /q} = c_{13}
\|y\|_{l_q^{\mathcal R_n}} = c_{13} \| \mathcal A^n f\|_{l_q^{\mathcal R_n}},
\end{equation*}
\begin{multline*} \tag{2.1.21}
\max_{j =1,\ldots,j_n}  \frac{ \| \mathcal E_{k +j} f\|_{L_p(I^d)} }
{ c_2 2^{-(k +j)} } 
 \ge c_{16} \left(\sum_{j =1}^{j_n} \left( \frac{ \|y_j\|_{l_p^{R_{k +j}} } }
{ 2^{-(k +j)(1 -(\alpha^{-1}, \e) /p +(\alpha^{-1}, \e) /q) +\epsilon j +\delta j}}
\right)^p \right)^{1/p} \\
 = c_{16} \left(\sum_{ \{\,(j,\nu,\lambda): j =1,\ldots,j_n,\ \nu \in
\Nu_{0, 2^{\kappa(k +j, \alpha)} -\e}^d,\ \lambda \in \Z_+^d(l -\e)\,\}}
\left| \frac {y_j^{\nu,\lambda}}
{\rho_j^{\nu,\lambda}} \right|^p \right)^{1 /p} \\
 = c_{16} \mu_{\mathcal B}(y) = c_{16} \mu_{\mathcal B}(\mathcal A^n f),
\end{multline*}
\begin{gather*}
\text{где } \mathcal B = \biggl\{\,y \in\R^{\mathcal R_n}:
\sum_{\{\,(j,\nu,\lambda): j =1,\ldots,j_n,\ \nu \in
\Nu_{0, 2^{\kappa(k +j, \alpha)} -\e}^d,\ \lambda \in \Z_+^d(l -\e)\,\}}
\left| \frac{ y_j^{\nu,\lambda}}
{\rho_j^{\nu,\lambda}} \right|^p \le 1 \,\biggr\}, \\
\text{а } \rho_j^{\nu,\lambda} = \rho_j = 2^{-(k +j)(1 -(\alpha^{-1}, \e) /p +
(\alpha^{-1}, \e) /q) +\epsilon j +\delta j}.
\end{gather*}

Объединяя (2.1.9), (2.1.20) и (2.1.21), в силу тех же соображений, что при
выводе (2.1.15), приходим к оценке
\begin{multline*} \tag{2.1.22}
b_n(C,X)
 \le 2 \sup_{ L \in \mathcal M_{R_k} \left(\mathfrak P_{k,k +j_n}^{\prime} \right) }
\inf_{ f \in L \setminus \{0\}} \frac{ c_{16} \| \mathcal A^n f\|_{l_q^{\mathcal R_n}} }
{c_{17} \mu_{\mathcal B}(\mathcal A^n f)} \\
 = c_{18} \sup_{ L \in \mathcal M_{R_k}\left(\mathfrak P_{k,k +j_n}^{\prime}\right) }
\inf_{ y \in \mathcal A^n (L) \setminus \{0\}} \frac{ \|y\|_{l_q^{\mathcal R_n}} }
{\mu_{\mathcal B}(y)} \\
 = c_{18} \sup_{\left\{\,M = \mathcal A^n(L): L \in \mathcal M_{R_k}
\left(\mathfrak P_{k,k +j_n}^{\prime}\right)\,\right \} }
\inf_{y \in M \setminus \{0\}} \frac{ \|y\|_{l_q^{\mathcal R_n}} }
{\mu_{\mathcal B}(y)} \\
 \le c_{18} \sup_{M \in \mathcal M_{R_k}\left(\R^{\mathcal R_n} \right) }
\inf_{y \in M \setminus \{0\}} \frac{ \|y\|_{l_q^{\mathcal R_n}} } {\mu_{\mathcal B}(y)} 
 = c_{18} b_{R_k}(\mathcal B, l_q^{\mathcal R_n}).
\end{multline*}

Беря $ \epsilon >0$ и $ \delta >0$ так, чтобы было
$ 1 -(\alpha^{-1}, \e) /p +(\alpha^{-1}, \e) /q -\epsilon -\delta >0$,
что можно сделать при условии $ 1 -(\alpha^{-1}, \e) /p +(\alpha^{-1}, \e) /q >0$,
и применяя (1.3.9), (1.2.5), находим, что
\begin{multline*} \tag{2.1.23}
b_{R_k}(\mathcal B, l_q^{\mathcal R_n})
 = \left(R_k \left(2^{-(k +1)(1 -(\alpha^{-1}, \e) /p +(\alpha^{-1}, \e) /q)
+\epsilon +\delta} \right)^{\frac{pq} {p -q}} \right)^{\frac{p -q} {pq}} \\
 \le \left(c_{19} 2^{k(\alpha^{-1}, \e)} \right)^{1/q -1/p} 2^{\epsilon +\delta}
\cdot 2^{-(k +1)(1 -(\alpha^{-1}, \e) /p +(\alpha^{-1}, \e) /q)} 
 \le c_{20} 2^{-k}.
\end{multline*}

Подставляя (2.1.23) в (2.1.22) и принимая во внимание (1.2.5), (2.1.7),
получаем (2.1.1) при $ p < q < \infty $ и условии
$ 1 -(\alpha^{-1}, \e) /p +(\alpha^{-1}, \e) /q >0$.

Проведенное доказательство сохраняет свою силу и при $ p < q = \infty$,
если соответствующим образом скорректировать (2.1.17) и (2.1.20), а
вместо (1.3.9) воспользоваться (1.3.9)$^\prime$.

Теперь в условиях теоремы 2.1.1 оценим $ b_n(B((H_p^\alpha)^\prime(Q)),L_q(Q)) $
через $ b_n(B((H_p^\alpha)^\prime(I^d)),L_q(I^d)), $ где
$ Q = x^0 +\delta I^d, \delta \in \R_+^d, x^0 \in \R^d. $
В этом случае для $ f \in B((H_p^\alpha)^\prime(Q)) $ имеем
\begin{equation*}
f = (h_{\delta,x^0})^{-1} h_{\delta,x^0} f = (h_{\delta,x^0})^{-1} \mathpzc f,
\end{equation*}
где в силу (1.4.3) функция $ \mathpzc f = h_{\delta,x^0} f \in c_{21}
B((H_p^\alpha)^\prime(I^d)), $ т.е.
\begin{equation*}
B((H_p^\alpha)^\prime(Q)) \subset (h_{\delta,x^0})^{-1} (c_{21}
B((H_p^\alpha)^\prime(I^d))).
\end{equation*}
Отсюда вследствие (1.3.2), (1.3.1), (1.3.3), (1.4.2),  при $ n \in \N $
выполняется неравенство
\begin{multline*} \tag{2.1.24}
b_n(B((H_p^\alpha)^\prime(Q)), L_q(Q)) \le b_n((h_{\delta,x^0})^{-1} (c_{21}
B((H_p^\alpha)^\prime(I^d))), L_q(Q)) = \\ c_{21} b_n((h_{\delta,x^0})^{-1}
(B((H_p^\alpha)^\prime(I^d))), L_q(Q)) \le \\
c_{21} \| (h_{\delta,x^0})^{-1} \|_{\mathcal B(L_q(I^d), L_q(Q))}
b_n(B((H_p^\alpha)^\prime(I^d)), L_q(I^d)) \le \\ 
c_{22} b_n(B((H_p^\alpha)^\prime(I^d)), L_q(I^d)).
\end{multline*}

Наконец, в условиях теоремы 2.1.1 установим оценку сверху величины
$ b_n(B((B_{p,\theta}^\alpha)^\prime(D)),L_q(D)). $
Для этого фиксируем $ \delta \in \R_+^d, x^0 \in \R^d $ такие, что $ D \subset
Q = x^0 +\delta I^d. $ Принимая во внимание (1.4.7), рассмотрим подпространство
\begin{equation*}
M = \{ (E^{d,\alpha,p,\infty,D} f) \mid_Q: f \in (H_p^\alpha)^\prime(D)\}
\subset (H_p^\alpha)^\prime(Q) \subset L_q(Q).
\end{equation*}
Тогда, благодаря (1.4.5) для $ f \in B((H_p^\alpha)^\prime(D)) $
имеют место соотношения
\begin{equation*}
f = (E^{d,\alpha,p,\infty,D} f) \mid_D = ((E^{d,\alpha,p,\infty,D} f) \mid_Q) \mid_D
= \mathpzc f \mid_D,
\end{equation*}
а $ \mathpzc f = (E^{d,\alpha,p,\infty,D} f) \mid_Q \in M $
ввиду (1.4.7) удовлетворяет условию
\begin{equation*}
\| \mathpzc f \|_{(H_p^\alpha)^\prime(Q)} \le
\| E^{d,\alpha,p,\infty,D} f \|_{(H_p^\alpha)^\prime(\R^d)} \le
c_{23} \| f \|_{(H_p^\alpha)^\prime(D)} \le c_{23},
\end{equation*}
т.е.
\begin{multline*}
B((H_p^\alpha)^\prime(D)) \subset U((c_{23} B((H_p^\alpha)^\prime(Q)))
\cap M) = U(c_{23} (B((H_p^\alpha)^\prime(Q)) \cap M)) = \\
c_{23}
U(B((H_p^\alpha)^\prime(Q)) \cap M),
\end{multline*}
где оператор $ U: M \cap L_q(Q) \mapsto (H_p^\alpha)^\prime(D) \cap L_q(D) \subset L_q(D) $
задаётся равенством $ U \mathpzc f = \mathpzc f \mid_D, \mathpzc f \in M. $
При этом, если для $ \mathpzc f \in M $ его образ $ U \mathpzc f = 0, $ то беря
$ f \in (H_p^\alpha)^\prime(D), $ для которого $ \mathpzc f =
(E^{d,\alpha,p,\infty,D} f) \mid_Q, $ ввиду (1.4.5) имеем
\begin{equation*}
f = (E^{d,\alpha,p,\infty,D} f) \mid_D = ((E^{d,\alpha,p,\infty,D} f) \mid_Q) \mid_D
= \mathpzc f \mid_D = U \mathpzc f = 0,
\end{equation*}
а, следовательно, $ \mathpzc f = (E^{d,\alpha,p,\infty,D} f) \mid_Q = 0 \mid_Q =0. $
Принимая во внимание эти обстоятельства, опираясь на (1.3.2), (1.3.1), (1.3.3),
при $ n \in \N $ получаем
\begin{multline*}
b_n(B((H_p^\alpha)^\prime(D)), L_q(D)) \le
b_n(c_{23} U(B((H_p^\alpha)^\prime(Q)) \cap M), L_q(D)) = \\
c_{23}
b_n(U(B((H_p^\alpha)^\prime(Q)) \cap M), L_q(D)) \le \\
c_{23}
\|U\|_{\mathcal B(M \cap L_q(Q), L_q(D)}
b_n(B((H_p^\alpha)^\prime(Q)) \cap M, M \cap L_q(Q)) \le \\
c_{23}
b_n(B((H_p^\alpha)^\prime(Q)) \cap M, M \cap L_q(Q)).
\end{multline*}
Из последнего неравенства, применяя (1.3.5), (1.3.2), выводим
\begin{multline*} \tag{2.1.25}
b_n(B((H_p^\alpha)^\prime(D)), L_q(D)) \le c_{23}
b_n(B((H_p^\alpha)^\prime(Q)) \cap M, M \cap L_q(Q)) = \\
c_{23}
b_n(B((H_p^\alpha)^\prime(Q)) \cap M,  L_q(Q)) \le c_{23}
b_n(B((H_p^\alpha)^\prime(Q)), L_q(Q)).
\end{multline*}

Учитывая, что в силу (1.1.2) имеет место включение
\begin{equation*}
B((B_{p,\theta}^\alpha)^\prime(D)) \subset c_{24}
B((H_p^\alpha)^\prime(D)),
\end{equation*}
используя (1.3.2), (1.3.1), (2.1.25), (2.1.24), находим, что
\begin{equation*}
b_n(B((B_{p,\theta}^\alpha)^\prime(D)), L_q(D)) \le c_{25}
b_n(B((H_p^\alpha)^\prime(I^d)), L_q(I^d)), n \in \N.
\end{equation*}

Соединение этой оценки с установленной выше оценкой правой части последнего
неравенства, завершает доказательство теоремы 2.1.1. $\square$
\bigskip

2.2. В этом пункте проводится оценка снизу $n$-поперечника по
Бернштейну шара $ B((B_{p,\theta}^\alpha)^\prime(D)) $ в пространстве
$L_q(D) $ для произвольной области $ D \subset \R^d. $ Для этого понадобятся следующие леммы.

Как видно из доказательства леммы 3.2.2 из [8] имеет место следующая
лемма.

Лемма 2.2.1

Пусть $ d \in \N $ и функция $ \phi \in C_0^\infty(I^d) $ такова,
что $ \phi \not\equiv 0 $. Для $ \kappa \in \N^d $ и множества
$  \Nu \subset \Nu_{0,\kappa -\e}^d $   обозначим через $ \L_{\kappa,\Nu} =
\L_{\kappa,\Nu}^{d,\phi} $ линейную оболочку системы функций
$ \{\phi(\kappa x -\nu), \nu \in \Nu\} $, а через $ \J_{\kappa,\Nu} =
\J_{\kappa,\Nu}^{d,\phi}: \R^{|\Nu|} \mapsto \L_{\kappa,\Nu}^{d,\phi} $ --
линейный оператор, определяемый равенством
$$
\J_{\kappa,\Nu} \beta = \sum_{\nu \in \Nu} \beta_\nu \phi(\kappa x -\nu),
\beta = \{\beta_\nu \in \R, \nu \in \Nu\} \in \R^{|\Nu|},
$$
а $ |\Nu| = \card \Nu $.

Тогда

1) для любого $ \kappa \in \N^d $ и любого непустого множества
$ \Nu \subset \Nu_{0,\kappa -\e}^d $ система функций $ \{\phi(\kappa x -\nu),
\nu \in \Nu\} $ --  линейно независима, и $ \J_{\kappa,\Nu} $ есть линейный
изоморфизм $ \R^{|\Nu|} $ на $ \L_{\kappa,\Nu} $;

2) для любого $ r: 1 \le r \le \infty $, существует  константа $ c_1(r,d,\phi) >0 $
такая, что для любых $ \kappa \in \N^d, \Nu \subset \Nu_{0,\kappa -\e}^d $
и $ \beta \in \R^{|\Nu|} $ справедливо равенство
\begin{equation*} \tag{2.2.1}
\| \J_{\kappa,\Nu} \beta \|_{L_r(\R^d)} = \| (\J_{\kappa,\Nu} \beta ) \mid_{I^d}\|_{L_r(I^d)} =
c_1 \kappa^{-r^{-1} \e} \|\beta\|_{l_r^{|\Nu|}};
\end{equation*}

3) для любых $ d \in \N, \alpha \in \R_+^d, \l = \l(\alpha), 1 \le p, \theta
< \infty $ существует константа $ c_2(d,\alpha,p,\theta,\phi) >0 $ такая, что
при любом $ k \in \N $ для $ \kappa = \Kappa(k,\alpha) \in \N^d $ с
компонентами $ \Kappa_j(k,\alpha) = [k^{1/\alpha_j}], j =1,\ldots,d $, для
любого $ \Nu \subset \Nu_{0,\kappa -\e}^d $ имеет место включение
\begin{equation*} \tag{2.2.2}
\J_{\kappa,\Nu}(\{ \beta \in \R^{|\Nu|}: \|\beta\|_{l_p^{|\Nu|}}
\le c_2 k^{-1} \kappa^{p^{-1} \e} \}) \subset
B((B_{p,\theta}^\alpha)^{\l}(\R^d)) \cap C_0^\infty(I^d).
\end{equation*}

Будем пользоваться следующим обозначением. Для множества $ S, $
состоящего из функций $ f, $ область определения которых содержит множество
$ D \subset \R^d, $ через $ S \mid_D $ обозначим множество $ S \mid_D =
\{ f \mid_D: f \in S\}. $

С помощью леммы 2.2.1 устанавливается лемма 2.2.2.

Лемма 2.2.2

Пусть $ d \in \N, \alpha \in \R_+^d, 1 \le p < \infty, 1 \le q  \le \infty $
удовлетворяют условию (1.2.10) и $ \theta \in \R: 1 \le \theta < \infty,
\l = \l(\alpha). $ Тогда существует константа $ c_3(d,\alpha,p,\theta,q) >0 $
такая, что для
$ C = (B((B_{p,\theta}^\alpha)^{\l}(\R^d)) \cap C_0^\infty(I^d)) \mid_{I^d},
X = L_q(I^d) $ для любого натурального числа $ n $ имеет место неравенство
\begin{equation*} \tag{2.2.3}
b_n(C,X) \ge c_3 n^{-1 /(\alpha^{-1}, \e) +p^{-1} -q^{-1}}
b_n(B(l_p^{2n}), l_q^{2n}).
\end{equation*}

Доказательство.

Фиксировав функцию $ \phi \in C_0^\infty(I^d): \phi \not\equiv 0, $
и выбрав для $ n \in \N $ число $ k \in \N: k \ge 2 $, так, чтобы для
$ \Kappa(k,\alpha) \in \N^d $ с   компонентами
$ \Kappa_j(k,\alpha) = [k^{1 /\alpha_j}], j =1,\ldots,d $, соблюдалось
условие
\begin{equation*} \tag{2.2.4}
(\Kappa(k -1,\alpha))^\e < 2n \le (\Kappa(k,\alpha))^\e,
\end{equation*}
для $ \kappa = \Kappa(k,\alpha) $ возьмём некоторое множество
$ \Nu \subset \Nu_{0,\kappa -\e}^d $, для которого $ |\Nu| = \card \Nu = 2n $,
и рассмотрим подпространство $ \L_{\kappa,\Nu}^{d,\phi} $, а также оператор
$ \J_{\kappa,\Nu}^{d,\phi} $ из леммы 2.2.1.

Тогда, используя сначала (2.2.2) и (1.3.2), а затем (1.3.1), (1.3.5),
и, далее, применяя (1.3.3) вместе с (2.2.1), и, наконец, принимая во
внимание, что для $ n $, удовлетворяющего (2.2.4), справедливо
соотношение $ n \asymp k^{(\alpha^{-1},
\e)} $, приходим к неравенству
\begin{multline*}
b_n(C,X) = b_n((B((B_{p,\theta}^\alpha)^{\l}(\R^d)) \cap C_0^\infty(I^d)) \mid_{I^d},
L_q(I^d)) \ge \\
b_n((\J_{\kappa,\Nu}(\{ \beta \in \R^{|\Nu|}: \|\beta\|_{l_p^{|\Nu|}}
\le c_2 k^{-1} \kappa^{p^{-1} \e} \})) \mid_{I^d}, L_q(I^d)) = \\
b_n((\J_{\kappa,\Nu}(c_2 k^{-1} \kappa^{p^{-1} \e} B(l_p^{|\Nu|})))
\mid_{I^d}, L_q(I^d)) = \\
c_2 k^{-1} \kappa^{p^{-1} \e} b_n((\J_{\kappa,\Nu}(B(l_p^{2n}))) \mid_{I^d},
L_q(I^d)) = \\
c_2 k^{-1} \kappa^{p^{-1} \e} b_n((\J_{\kappa,\Nu}(B(l_p^{2n}))) \mid_{I^d},
\L_{\kappa,\Nu} \mid_{I^d} \cap L_q(I^d)) \ge \\
c_4 k^{-1} \kappa^{p^{-1} \e} \kappa^{-q^{-1} \e} b_n(B(l_p^{2n}), l_q^{2n}) \ge 
c_3 n^{-1 /(\alpha^{-1}, \e) +p^{-1} -q^{-1}}
b_n(B(l_p^{2n}), l_q^{2n}). \square
\end{multline*}

Теорема 2.2.3

Пусть $ d \in \N, \alpha \in \R_+^d, D $ -- ограниченная область $ \alpha $-типа
в $ \R^d, 1 \le p < \infty, 1 \le q \le \infty $ удовлетворяют условию (1.2.10)
и $ \theta: 1 \le \theta \le \infty. $ Тогда существует константа
$ c_5(d,\alpha,D,p,\theta,q) >0 $ такая, что для $ C =
B((B_{p,\theta}^\alpha)^\prime(D)), X = L_q(D) $ для любого натурального числа
$ n $ выполняется неравенство
\begin{multline*} \tag{2.2.5}
b_n(C,X) \ge c_5 \cdot
\begin{cases}
n^{-1 /(\alpha^{-1}, \e)}, \parbox[t]{7cm}{ при $ 1 \le q \le p \le 2$ или $ 1 \le q = p <
\infty$ или  $ 1 \le p < q \le \infty$ ;} \\
n^{-\!(1 /(\alpha^{-1}, \e) -1 /p +1/2)},  \text{ при $ 1 \le q \le 2 < p < \infty $;} \\
n^{-\!(1 /(\alpha^{-1}, \e) -1 /p +1 /q)},  \text{ при $ 2 \le q < p < \infty.$}
\end{cases}
\end{multline*}

Доказательство.

Ввиду (1.1.2), (1.3.2) соотношение (2.2.5) достаточно проверить при $ \theta \ne \infty. $

Фиксируем точку $ x^0 \in \R^d $ и вектор $ \delta \in \R_+^d $
такие, что $ Q = (x^0 +\delta I^d) \subset D. $
Сначала заметим, что в условиях леммы при $ n \in \N, \l = \l(\alpha), $
благодаря (2.2.3), (1.1.3) и (1.3.2), имеет место неравенство
\begin{multline*} \tag{2.2.6}
c_3 n^{-1 /(\alpha^{-1}, \e) +(p^{-1} -q^{-1}) } b_n(B(l_p^{2n}), l_q^{2n}) \le \\
b_n((B((B_{p,\theta}^\alpha)^\prime(\R^d)) \cap C_0^\infty(I^d)) \mid_{I^d}, L_q(I^d)).
\end{multline*}

Далее, принимая во внимание, что в силу (1.4.4) справедливо включение
\begin{multline*}
(B((B_{p,\theta}^\alpha)^\prime(\R^d)) \cap C_0^\infty(I^d)) \mid_{I^d}
\subset \\
h_{\delta,x^0} (\{ (F \mid_Q):
F \in (c_6 B((B_{p,\theta}^\alpha)^\prime(\R^d))) \cap C_0^\infty(Q) \}) = \\
c_6 h_{\delta,x^0} (\{ (F \mid_Q):
F \in B((B_{p,\theta}^\alpha)^\prime(\R^d)) \cap C_0^\infty(Q) \}),
\end{multline*}
благодаря (1.3.2), (1.3.1), (1.3.3), (1.4.1), имеем
\begin{multline*} \tag{2.2.7}
b_n((B((B_{p,\theta}^\alpha)^\prime(\R^d)) \cap C_0^\infty(I^d)) \mid_{I^d}, L_q(I^d)) \le \\
b_n(c_6 h_{\delta,x^0} (\{ (F \mid_Q):
F \in B((B_{p,\theta}^\alpha)^\prime(\R^d)) \cap C_0^\infty(Q) \}),L_q(I^d)) \le \\
c_6 \delta^{-q^{-1} \e} b_n(\{ F \mid_Q:
F \in B((B_{p,\theta}^\alpha)^\prime(\R^d)) \cap C_0^\infty(Q) \},L_q(Q)).
\end{multline*}

Теперь, определяя в $ L_q(D) $ замкнутое подпространство
$$
L_q^0(Q,D) = \{f \in L_q(D): f = \chi_Q f\},
$$
заметим, что
$$
\{ (F \mid_D):
F \in B((B_{p,\theta}^\alpha)^\prime(\R^d)) \cap C_0^\infty(Q) \} \subset
L_q^0(Q,D),
$$
а для оператора $ U: L_q(D) \ni f \mapsto U f = f \mid_Q \in L_q(Q), $
его сужение $ U \mid_{L_q^0(Q,D)} $ является изометрическим изоморфизмом
пространства $ L_q^0(Q,D) \cap L_q(D) $ на $ L_q(Q). $
С учётом сказанного применяя (1.3.3), а затем (1.3.5) и (1.3.2), получаем
\begin{multline*} \tag{2.2.8}
b_n(\{ F \mid_Q:
F \in B((B_{p,\theta}^\alpha)^\prime(\R^d)) \cap C_0^\infty(Q) \},L_q(Q)) = \\
b_n(U \mid_{L_q^0(Q,D)} (\{ (F \mid_D):
F \in B((B_{p,\theta}^\alpha)^\prime(\R^d)) \cap C_0^\infty(Q) \}),L_q(Q)) \le \\
b_n(\{ (F \mid_D):
F \in B((B_{p,\theta}^\alpha)^\prime(\R^d)) \cap C_0^\infty(Q) \}, L_q^0(Q,D) \cap L_q(D)) = \\
b_n(\{ (F \mid_D):
F \in B((B_{p,\theta}^\alpha)^\prime(\R^d)) \cap C_0^\infty(Q) \}, L_q(D)) \le \\
b_n(\{ f: f \in B((B_{p,\theta}^\alpha)^\prime(D)),L_q(D)).
\end{multline*}

Соединяя (2.2.8), (2.2.7), (2.2.6), видим, что
\begin{equation*} \tag{2.2.9}
b_n(B((B_{p,\theta}^\alpha)^\prime(D)),L_q(D)) \ge
c_7 n^{-1 /(\alpha^{-1}, \e) +(p^{-1} -q^{-1}) } b_n(B(l_p^{2n}), l_q^{2n}).
\end{equation*}

Подставляя в (2.2.9) оценку (1.3.10), приходим к (2.2.5). $\square$

\end{document}